\documentclass[11pt,reqno]{amsart}
\usepackage{graphicx,enumerate,amssymb,bbold}
\usepackage[notrig]{physics}
\usepackage{amsmath}
\usepackage{tikz-cd}
\usepackage{mathrsfs}
\usepackage[font=footnotesize,labelfont=small]{subcaption}
\usepackage[ruled,vlined]{algorithm2e}
\usepackage[left=3.5cm,right=3.5cm,top=2.5cm,bottom=2.5cm,includeheadfoot]{geometry} 
\usepackage[colorlinks=true, pdfstartview=FitV, linkcolor=blue, 
            citecolor=blue, urlcolor=blue]{hyperref}

\usepackage[disable]{todonotes}

\makeatletter
\providecommand\@dotsep{5}
\def\listtodoname{List of Todos}
\def\listoftodos{\@starttoc{tdo}\listtodoname}
\makeatother

\numberwithin{equation}{section}
\numberwithin{figure}{section}
\theoremstyle{plain}
\newtheorem{theorem}{Theorem}[section]
\newtheorem{lemma}[theorem]{Lemma}
\newtheorem{proposition}[theorem]{Proposition}
\newtheorem{corollary}[theorem]{Corollary}

\theoremstyle{definition}
\newtheorem{definition}[theorem]{Definition}

\newtheorem{remark}[theorem]{Remark}

\newcommand{\om}{\omega}
\newcommand{\what}{\widehat}
\newcommand{\wtilde}{\widetilde}

\newcommand{\vvar}{\mrm{Var}}

\newcommand{\ggradg}{\ggrad^{g}\!}

\newcommand{\bitem}{\begin{itemize}}
\newcommand{\eitem}{\end{itemize}}
\newcommand{\mc}[1]{\mathcal{#1}}

\newcommand{\F}{\mathbb{F}}
\newcommand{\N}{\mathbb{N}}
\newcommand{\R}{\mathbb{R}}

\newcommand{\EE}{\mathbb{E}}

\newcommand{\bpm}{\begin{pmatrix}}
\newcommand{\epm}{\end{pmatrix}}
\newcommand{\bvm}{\begin{vmatrix}}
\newcommand{\evm}{\end{vmatrix}}
\newcommand{\bsm}{\left(\begin{smallmatrix}}
\newcommand{\esm}{\end{smallmatrix}\right)}
\newcommand{\T}{\top}

\newcommand{\la}{\langle}
\newcommand{\ra}{\rangle}

\newcommand{\mrm}[1]{\mathrm{#1}}

\newcommand{\veps}{\varepsilon}

\newcommand{\w}{\omega}
\newcommand{\Om}{\Omega}

\newcommand{\vphi}{\varphi}

\newcommand{\eins}{\mathbb{1}}

\DeclareMathSymbol{\mydiv}{\mathbin}{symbols}{"04}

\DeclareMathOperator{\Diag}{Diag}

\DeclareMathOperator{\rint}{rint}

\DeclareMathOperator{\argmax}{arg max}

\DeclareMathOperator{\ggrad}{grad}

\DeclareMathOperator{\Exp}{Exp}

\newcommand{\Mcrit}{{\mc{M}_\mrm{crit}}}
\newcommand{\Mreg}{{\mc{Q}}}

\newcommand{\SM}{{\mc{S}}}
\newcommand{\WM}{{\mc{W}}}

\newcommand{\BS}{{\eins_\SM}}
\newcommand{\BW}{{\eins_\WM}}

\newcommand{\TS}{{T_{0}}}

\newcommand{\TW}{{\mc{T}_{0}}}

\newcommand{\PTS}{{P_\TS}}
\newcommand{\PTW}{{\mc{P}_\TW}}

\newcommand{\ROS}{{R}}
\newcommand{\ROW}{{\mc{R}}}

\newcommand{\Graph}{{\mc{G}}}
\newcommand{\Nodes}{{\mc{V}}}
\newcommand{\Edges}{{\mc{E}}}

\newcommand{\nhood}{{\mc{N}}}

\title{On the Geometric Mechanics of Assignment Flows \\ for Metric Data Labeling}
\author[F.~Savarino, P.~Albers, C.~Schn\"{o}rr]{Fabrizio Savarino, Peter Albers, Christoph Schn\"{o}rr}
\address[F.~Savarino]{STRUCTURES Excellence Cluster, Heidelberg University, Germany}
\urladdr{\url{https://www.thphys.uni-heidelberg.de/~structures/}}
\address[P.~Albers]{Mathematical Institute, Heidelberg University, Germany}
\urladdr{\url{https://www.mathi.uni-heidelberg.de/~palbers/}}
\address[C.~Schn\"{o}rr]{Institute of Applied Mathematics, Heidelberg University, Germany} 
\urladdr{\url{https://ipa.math.uni-heidelberg.de}}
\date{} 


\begin{document}

\begin{abstract}
Metric data labeling refers to the task of assigning one of multiple predefined labels to every given datapoint based on the metric distance between label and data. This assignment of labels typically takes place in a spatial or spatio-temporal context. Assignment flows are a class of dynamical models for metric data labeling that evolve on a basic statistical manifold, the so called assignment manifold, governed by a system of coupled replicator equations. In this paper we generalize the result of a recent paper for uncoupled replicator equations and adopting the viewpoint of geometric mechanics, relate assignment flows to critical points of an action functional via the associated Euler-Lagrange equation. We also show that not every assignment flow is a critical point and characterize precisely the class of coupled replicator equations fulfilling this relation, a condition that has been missing in recent related work. Finally, some consequences of this connection to Lagrangian mechanics are investigated including the fact that assignment flows are, up to initial conditions of measure zero, reparametrized geodesics of the so-called Jacobi metric.

\end{abstract}

\maketitle
\tableofcontents

\section{Introduction}
\label{sec:Introduction}

\subsection{Overview, Motivation}\label{sec:Overview}
\textit{Semantic image segmentation}, a.k.a.~\textit{image labeling}, denotes the problem to partition an image into meaningful parts. Applications are abound and include interpretation of traffic scenes by computer vision systems, medical image analysis, remote sensing, etc. The state of the art is based on deep networks that were trained on very large data sets. A recent survey \cite{Minaee:2021tx} reviews a vast number of different deep network architectures and their empirical performance on various benchmark data sets. Among the challenges discussed in \cite[Sec.~6.3]{Minaee:2021tx}, the authors write: ``... a concrete study of the underlying behavior / dynamics of these models is lacking. A better understanding of the theoretical aspects of these models can enable the development of better models curated toward various segmentation scenarios.''

In \cite{Astroem2017}, a class of dynamical systems for image labeling, called \textit{assignment flows}, was introduced in order to contribute to the mathematics of deep networks and learning. We refer to Section \ref{sec:Mechanics_of_Assignment_Flows} for a precise definition and to \cite{Schnorr2019aa} for a review of recent related work. \textit{Assignment flows} correspond to solutions $W(t)$ of a high-dimensional system of coupled \textit{ordinary differential equations} (ODEs) of the form
\begin{equation}\label{eq:AF}
\dot W(t) = \mc{R}_{W(t)}[F\big(W(t)\big)], 
\end{equation}
that evolve on the so-called \textit{assignment manifold} $\mc{W}$. Each ODE of this system is a \textit{replicator equation} \cite{hofbauer2003,sandholm2010}
\begin{equation}\label{eq:replicator}
\dot W_{i} = R_{W_{i}}F_{i}(W),\qquad
\dot W_{ij} = W_{ij}\Big(F_{j}(W)-\sum_{l=1}^cW_{il}F_{l}(W)\Big),\quad j\in\{1,\ldots,n\},
\end{equation}
whose solution $W_{i}(t)\in\mc{S}:=\rint\Delta_{n-1}\subset \R_{+}^{n}$ evolves on the relative interior of the probability simplex that is equipped with the Fisher-Rao metric $g$ \cite{amari2007methods} and is labeled by a vertex $i\in\mc{V}$ of an underlying graph $\mc{G}=(\mc{V},\mc{E})$. The assignment manifold $\mc{W}=\mc{S}\times\dotsb\times\mc{S}$ is the product of the Riemannian manifolds $(\mc{S},g)$ with respect to all vertices $i\in\mc{V}$.

The essential component of the vector field of \eqref{eq:replicator} are a collection of \textit{affinity functions} $F_{ij}\colon \mc{W}\to\R$ that measure the affinity (fitness, etc.) of the classes (types, species, etc.) $j\in[n]$. The differences of these affinity values to their expected (or average) value on the right-hand side of \eqref{eq:replicator}, together with the multiplication by $W_{ij}$, define the replicator equation. For suitably defined affinity functions, the solution of this equation is supposed to perform a \textit{selection} of some class $j$: $W_{i}(t)$ converges for $t\to\infty$ to a vertex of $e_{j} \in \Delta_{n-1}$ and in this sense encodes the \textit{decision} to \textit{assign} the class label $j$ to the vertex $i\in\mc{V}$ and to any data indexed by $i$, like e.g.~the color value in some image, see Section \ref{sec:AF} for more details.

The basic idea underlying the \textit{assignment flow approach} \eqref{eq:AF} is to assign a replicator equations to \textit{each} vertex of an underlying  graph and to \textit{couple} them through smooth nonlinear  interactions of the assignment vectors $\{W_{k}\colon k\in\mc{N}_{i}\subset\mc{V}\}$ within neighborhoods $\mc{N}_{i}$ around each vertex $i\in\mc{V}$. This is why the argument of $F_{i}$ in \eqref{eq:replicator} is $W$ rather than $W_{i}$. As a consequence, dynamic label assignments are performed by solving \eqref{eq:AF} at each vertex $i$ depending on the \textit{context} in terms of all other decisions. The fact that $W(t)$ assigns class labels at each vertex when $t\to\infty$ is not clear \textit{a priori} but depends on $F$. We refer to \cite{Zern:2020aa} for the study of a basic instance of $F$ and sufficient conditions that ensure unique labeling decisions.

The \textit{connection to deep networks} results from approximating the flow by geometric integration. The simplest such scheme among a range of proper schemes \cite{Zeilmann:2020aa}, the geometric Euler scheme with discrete time index $t$ and stepsize $h_{(t)}$, yields the iterative update rule
\begin{equation}\label{eq:AF-Euler}
W_{i}^{(t+1)} = \Exp_{W_{i}^{(t)}} \circ \ROS_{W_{i}^{(t)}}\big(h_{(t)} F_{i}(W^{(t)})\big),\quad i\in\mc{V},
\end{equation}
where $\Exp\colon T\mc{W}\to\mc{W}$ denotes the exponential map of the so-called e-connection of information geometry \cite{amari2007methods,ay2017information}. The key observation to be made here is that for the choice of a linear affinity map $F$ the right-hand side of \eqref{eq:AF-Euler} involves the two essential ingredients of most deep network architectures: 
\begin{enumerate}
\item
A linear operation at each vertex of the underlying graph parametrized by network parameters, here given as part of the definition of the linear affinity map $F$. 
\item 
A pointwise smooth nonlinearity, here given by the exponential and replicator maps $\Exp_{W_{i}} \circ \ROS_{W_{i}}$.
\end{enumerate}

The connection between general continuous-time ODEs and deep networks has been picked out as a central them by \cite{Haber:2017aa,Chen:2018ab} and classifies the assignment flow as a particular `Neural ODE'. The above-mentioned limited understanding of what deep networks really do underlines the importance of characterizing and understanding the dynamics \eqref{eq:AF} of assignment flows.

\subsection{Contribution, Organization}

The aim of this paper is to exhibit a natural Lagrangian $L\colon T\mc{W}\to\R$ of the form \textit{kinetic minus potential energy} and to characterize solutions $W(t)$ to \eqref{eq:AF} as stationary points of a corresponding action functional
\begin{equation}\label{eq:def-mcL}
\mc{L}(W) = \int_{0}^{t}L(W,\dot W) dt.
\end{equation}
Our result generalizes the result of a recent paper \cite{raju2018variational}, where an action functional was introduced for the evolution $p(t)$ of a \textit{single} discrete probability vector on the corresponding probability simplex. By contrast, equation \eqref{eq:AF} couples the evolution of a (typically large) number of assignment vectors across the underlying graph. In particular, we characterize precisely the admissible class of affinity functions $F$ that establishes the connection between \eqref{eq:def-mcL} and the corresponding Euler-Lagrange equation, a condition that is missing in \cite{raju2018variational}, see also Section \ref{sec:RK18}. Furthermore, we show that except for starting points in a specific set of measure zero, the set of \textit{Ma\~n\'e critical points}, solutions of the assignment flow are reparametrized geodesics of the so called \textit{Jacobi metric}. Finally, using the Legendre transform, we compute an explicit expression of the Hamiltonian system associated to \eqref{eq:def-mcL} in the form of the equivalent Lagrangian system on $T\WM$.

The paper is organized as follows. Section \ref{sec:Lagrangian} collects basic notions of geometric mechanics that are required in the remainder of the paper. The assignment flow and our novel results are presented in Section \ref{sec:Mechanics_of_Assignment_Flows}, followed by a discussion in Section \ref{sec:Discussion}. We conclude in Section \ref{sec:Conclusion}.

\subsection{Basic Notation}
In accordance with the standard notation in differential geometry, coordinates of vectors have upper indices. For any $k \in \N$, we set $[n] := \{ 1, \ldots, n\} \subset \N$. The standard basis of $\R^d$ is denoted by $\{e_1, \ldots, e_d\}$ and we set $\eins_d := (1, \ldots, 1)^\T \in \R^d$. 

Depending on the arguments, $\la a, b\ra$ denotes the Euclidean inner product of vectors or the inner product $\la A, B\ra=\tr(A^{\T} B)$ of matrices inducing the Frobenius norm $\|A\|_{F}=\la A, A\ra^{1/2}$. 
The identity matrix is denoted by $I_d \in \R^{d\times d}$ and the $i$-th row vector of any matrix $A$ by $A_i$.

The \textit{linear} dependence of a mapping $F$ on its argument $x$ is indicated by square brackets $F[x]$, if $F$ is just a matrix we simply write $Fx$. The \textit{adjoint} of a linear operator $F\colon \R^{m\times n} \to \R^{m\times n}$ with respect to the standard matrix inner product on $\R^{m\times n}$ is denoted by $F^*$ and fulfills
\begin{equation}\label{eq:def-adjoint-lin-op}
  \la F^*[A], B\ra = \la A, F[B]\ra, \quad \text{for all } A, B \in \R^{m\times n}
\end{equation}

Inequalities between vectors or matrices are to be understood componentwise. For $a, b \in \R^d$, we denote componentwise multiplication (Hadamard product) by
\begin{equation}
  a \diamond b := (a^1b^1, \ldots, a^d b^d)^\T
\end{equation}
and componentwise division, for $b\neq 0$, by $\frac{a}{b} = (\frac{a^1}{b^1}, \ldots, \frac{a^d}{b^d})^\T$.  We further set 
\begin{equation}
  a^{\diamond k} := a^{\diamond (k-1)}\diamond a\quad \text{and}\quad a^{\diamond 0} := \eins_d.
\end{equation}
Finally, if $p \in \R^d$ is a probability vector, i.e.~$p \geq 0$ and $\la p, \eins_d\ra = 1$, then the expected value and variance of a vector $a \in \R^d$ (interpreted as a random variable $a \colon [d] \to \R$) is 
\begin{equation}
  \EE_p[a] = \la p, a\ra \quad \text{and} \quad \vvar_p(a) = \EE_p[a^2] - (\EE_p[a])^2 = \la p, a^{\diamond 2}\ra - \la p, a\ra^2.
\end{equation}

\section{Elements from Geometric Mechanics}
\label{sec:Lagrangian}

In this section, we collect from \cite[Ch.~3]{abraham1987foundations} some basic notions of geometric mechanics that are required in subsequent sections.

\subsection{Hamiltonian Systems.}\label{sec:Lagrangian:Hamiltonian-Systems}
Let $(N, \om)$ be a symplectic manifold with the symplectic two-form $\w$, and let $H \colon N \to \R$ be a smooth function, called the \textit{Hamiltonian}. The \textit{Hamiltonian vector field $X_H$} corresponding to $H$ is defined as symplectic gradient by
\begin{equation}\label{eq:dH-om}
  dH|_x [v] = \om_x(X_H(x), v), \quad \text{for all } x\in N, v \in T_xN.
\end{equation}
The triplet $(N, \om, X_H)$ is called a \textit{Hamiltonian system}. By \cite[Prop.~3.3.2]{abraham1987foundations}, a curve $\gamma(t)$ is an integral curve of $X_H$, i.e.
\begin{equation}
  \dot{\gamma(t)} = X_H(\gamma(t)),
\end{equation}
if and only if in Darboux coordinates $(q^1, \ldots, q^n, p_1, \ldots, p_n)$ for $\om$, the \textit{Hamiltonian equations} hold for the curve $\gamma(t) = (q(t), p(t))$,  
\begin{equation}
  \dot{q}^i(t) = \frac{\partial H}{\partial p_i} (q(t), p(t))\quad \text{and} \quad \dot{p}_i = - \frac{\partial H}{\partial q^i} (q(t), p(t)), \quad \text{for all } i \in [n].
\end{equation}
The value of the Hamiltonian $H(\gamma(t))$ (also called \textit{energy}) is constant along integral curves of $X_H$.

For any smooth manifold $M$, the cotangent bundle $(T^*M,\om^\mrm{can})$ is a basic instance of the above situation, with the canonical symplectic form $\om^\mrm{can}$.
Thus any smooth function $H \colon T^*M \to \R$ gives rise to a Hamiltonian system, where $T^*M$ is interpreted as \textit{momentum phase space} and $H$ represents an energy.

\subsection{Lagrangian Systems.} \label{sec:Lagrangian:Lagrangian-Systems}
Suppose $M$ is a smooth manifold. Similar to Hamiltonian systems on momentum phase space $T^*M$, there is a related concept on the tangent bundle $TM$, interpreted as \textit{velocity phase space}. In this context,  a smooth function $L \colon TM \to \R$ is called \textit{Lagrangian}. For a given point $x \in M$, denote the restriction of $L$ to the fiber $T_xM$ by $L_x := L|_{T_xM} \colon T_xM \to \R$. The \textit{fiber derivative} of $L$ is defined as
\begin{equation}\label{eq:def-Fiber-derivative}
  \F L \colon TM \to T^*M, \quad (x, v) \mapsto \F L(x,v) := (x, dL_x|_v),
\end{equation}
where $dL_x|_v \colon T_xM \to \R$ is the differential of $L_x$ at $v \in T_xM$. The function $L$ is called a \textit{regular Lagrangian} if $\F L$ is regular at all points (i.e.~$\F L$ is a submersion), which is equivalent to $\F L\colon TM \to T^*M$ being a local diffeomorphism \cite[Prop.~3.5.9]{abraham1987foundations}. Furthermore, $L$ is called \textit{hyperregular Lagrangian} if $\F L \colon TM \to T^*M$ is a diffeomorphism. A class of hyperregular Lagrangians that will be relevant in Section \ref{sec:Mechanics_of_Assignment_Flows}, is given as equation \eqref{eq:Lagrangian-Systems:L-Riemann-mfld} below.

The \textit{Lagrangian two-form} $\om_L$ is defined as the pullback
of the canonical symplectic form $\om^\mrm{can}$ on the cotangent bundle $T^*M$ under the fiber derivative $\F L$
\begin{equation}\label{eq:Lagrangian-two-form-general-def}
  \om_L := (\F L)^*\om^\mrm{can}.
\end{equation}
According to \cite[Prop.~3.5.9]{abraham1987foundations}, $\om_L$ is a symplectic form on $T^*M$ if and only if $L$ is a regular Lagrangian. In the following, we only consider regular Lagrangians. The \textit{action} associated to a Lagrangian $L \colon TM \to \R$ is defined by
\begin{equation}
  A \colon TM \to \R, \quad (x, v) \mapsto \F L(x, v)[v] = dL_x|_v [v],
\end{equation}
and the \textit{energy function} by $E:= A - L$, that is
\begin{equation}\label{eq:def-energy}
  E\colon TM \to \R,\quad (x, v) \mapsto \F L(x, v)[v] - L(x, v) = dL_x|_v[v] - L(x, v).
\end{equation}
The \textit{Lagrangian vector field for $L$} is the unique vector field $X_E$ on $TM$ satisfying
\begin{equation}\label{eq:def-Lagrangian-vector-field}
  dE|_x[v] = \om_{L, x}(X_E, v) \quad \text{for all } x \in M, v \in T_x M.
\end{equation}
Since we assume $L$ to be regular, $X_E$ is nothing else than the symplectic gradient of $L$ with respect to $\om_L$.
A curve $\gamma(t) = (x(t), v(t))$ on $TM$ is an integral curve of $X_E$, i.e.
\begin{equation}\label{eq:X_E-integral-curve}
  \dot{\gamma}(t) = X_E(\gamma(t)),
\end{equation}
if $v(t) = \dot{x}(t)$ and the classical Euler-Lagrange equations in local coordinates 
\begin{equation}\label{eq:EL-classical}
  \frac{d}{dt} \bigg( \frac{\partial L}{\partial \dot x^i} \big(x(t), \dot{x}(t)\big)\bigg) = \frac{\partial L}{\partial x^i}\big(x(t), \dot{x}(t)\big) \quad \text{for all } i \in [n]
\end{equation}
are satisfied. 
Let $\gamma \colon I \to TM$ be any integral curve of $X_E$. Then 
\begin{equation}\label{eq:E-constant-general}
\tfrac{d}{dt}E(\gamma)=0, 
\end{equation}
that is the energy $E$ is constant along $\gamma$, analogous to the constancy of the Hamiltonian $H$ due to \eqref{eq:dH-om}. The subsequent Section \ref{sec:Legendre-transform} makes this connection explicit.

\subsection{The Legendre Transform.}\label{sec:Legendre-transform}
Let $L\colon TM\to\R$ be a hyperregular Lagrangian, i.e.~the fiber derivative $\F L \colon TM \to T^*M$ is a diffeomorphism. Then the Lagrangian system on $TM$ and the Hamiltonian system on $T^*M$ are related to each other by the \textit{Legendre transformation}, with the Hamiltonian $H\colon T^{\ast}M\to\R$ corresponding to the energy $E$ via
\begin{equation}\label{eq:H-by-E}
  H = E \circ (\F L)^{-1}.
\end{equation}
Accordingly, the Hamiltonian vector field $X_H$ on $T^*M$ and the Lagrangian vector field $X_E$ on $TM$ are $\F L$ related \cite[Thm.~3.6.2]{abraham1987foundations}, that is
\begin{equation}
  (\F L)_* X_E = X_H,
\end{equation}
and thus integral curves of $X_E$ are mapped to integral curves of $X_H$ and vice versa. Furthermore, the base integral curves of $X_E$ and $X_H$ coincide.

Therefore, as a consequence of \eqref{eq:H-by-E} and for a hyperregular Lagrangian $L$, the energy $E$ is just another representation of the corresponding Hamiltonian $H$.

\subsection{Mechanics on Riemannian Manifolds.}\label{sec:Mechanics-Riemannian}
Let $(M, h)$ be a Riemannian manifold. Suppose a smooth function $G \colon M \to \R$, called \textit{potential}, is given and consider the Lagrangian 
\begin{equation}\label{eq:Lagrangian-Systems:L-Riemann-mfld}
  L(x, v) = \tfrac{1}{2}\|v\|_h^2 - G(x), \quad (x, v) \in TM.
\end{equation}
It then follows (see \cite[Sec.~3.7]{abraham1987foundations} or by direct computation) that the fiber derivative of $L$ is the canonical isomorphism
\begin{equation}
  \F L = h^\flat \colon TM \to T^*M.
\end{equation}
Hence the Lagrangian $L$ is hyperregular with action $A$ and energy $E = A - L$ given by
\begin{equation}\label{eq:Lagrangian-Systems:L-Riemann-mfld-Action-and-Energy}
  A(x, v) = \|v\|_h^2 \quad  \text{and} \quad E(x, v) = \tfrac{1}{2}\|v\|_h^2 + G(x) \quad \text{for all } (x, v) \in TM.
\end{equation}

\begin{proposition}\label{prop:Lagrangian-Systems:E-L-Eq-on-Rmflds}
  \textbf{(\cite[Prop.~3.7.4]{abraham1987foundations}).}
  Let $(M, h)$ be a Riemannian manifold, $\pi \colon TM \to M$ the natural projection and $L\colon TM\to\R$ the Lagrangian defined by \eqref{eq:Lagrangian-Systems:L-Riemann-mfld}. Then the curve $\gamma \colon I \to TM$ with $\gamma(t) = (x(t), v(t))$ is an integral curve of the Lagrangian vector field $X_E$, i.e.~satisfies the Euler-Lagrange equation, if and only if the corresponding base integral curve $\pi\circ \gamma = x \colon I \to M$ satisfies
  \begin{equation}
    D^h_t\dot{x}(t) = - \ggrad^h G(x(t)),
  \end{equation}
\end{proposition}
\noindent
where $D^h_t = \nabla^{h}_{\dot x}$ is the covariant derivative along $x$ with respect to the Riemannian (Levi-Civita) connection $\nabla^{h}$.
Here, $\mrm{grad}^{h}G$ denotes the Riemannian gradient of the potential $G$.

\section{Mechanics of Assignment Flows}
\label{sec:Mechanics_of_Assignment_Flows}

In this section, we get back to the scenario of image labeling informally introduced in Section~\ref{sec:Overview}. Section~\ref{sec:AF} completes the definition of the assignment flow approach \eqref{eq:AF}. The assignment manifold underlying the assignment flow is introduced in Section~\ref{sec:AssManifold} together with Fisher-Rao metric in Section~\ref{sec:FR-Metric}. We state and prove the main result of this paper in Section~\ref{sec:action-functional}.

\subsection{Assignment Manifold} \label{sec:AssManifold}

Let $\Graph = (\Nodes, \Edges)$ denote an undirected graph and identify
\begin{equation}
  \Nodes = [m]\quad \text{with} \quad m := |\Nodes|.
\end{equation}
Assume that for every node $i \in \Nodes$ some data point $f_i$ is given in a metric space $(\mc{F}, d_{\mc{F}})$, together with a set $\mc{F}_* = \{f_1^*, \ldots, f_n^*\} \subset \mc{F}$ of predefined prototypes, also called \textit{labels}, identified with 
\begin{equation}
  \mc{F}_* = [n] \quad \text{for}\quad n := |\mc{F}_*|.
\end{equation}
Context based \textit{metric data classification} or \textit{labeling} refers to the task of assigning to each node $i \in \Nodes$ a suitable label in $\mc{F}_*$, based on the metric distance to the given data $f_i$ and the relation between data points encoded by the edge set $\Edges$.

As introduced in Section~\ref{sec:Overview}, for every $i \in \Nodes$ the \textit{assignment} of labels $\mc{F}_*$ to a data point $f_{i}$ is represented by an assignment vector $W_{i}(t)$, where the $j$-th entry $W_i^j(t)$ represents the probability for the $j$-th label $f_{j}^{\ast}$. These assignment vectors are determined by \eqref{eq:replicator} and evolve on the relative interior of the $(n-1)$-simplex
\begin{equation}\label{eq:Introduction:def-S}
  \SM := \{ p \in \R^n\colon p > 0 \text{ and } \la p, \eins_n\ra = 1\}
\end{equation}
with barycenter
\begin{equation}
\eins_{\mc{S}} := \tfrac{1}{n}\eins_{n}.
\end{equation}
Accordingly, all probabilistic label choices on the graph are encoded as a single point $W\in\mc{W}$ on the product space
\begin{equation}\label{eq:Introduction:def-W}
  \WM := \SM \times \ldots \times \SM\qquad \text{($m = |\Nodes|$ factors)}, 
\end{equation}
with barycenter
\begin{equation}\label{eq:def-barycenter-W-product}
  \BW := (\eins_{\mc{S}},\dotsc,\eins_{\mc{S}})^{\T}.
\end{equation}
Thus, the $i$-th component of $W = (W_k)_{k\in \Nodes}$ represents the probability distribution of label assignments for node $i \in \Nodes$
\begin{equation}\label{eq:component_i_of_W}
  W_i=(W_{i}^1,\dotsc,W_{i}^n)^{\T} \in \SM.
\end{equation}
In the following, we always identify the space $\WM$ from \eqref{eq:Introduction:def-W} with its matrix embedding
\begin{equation}\label{eq:Mechanics_of_Assignment_Flows:embedding-W}
  \WM = \{ W \in \R^{m\times n}\colon W > 0 \text{ and } W\eins_n = \eins_m\},
\end{equation}
by regarding the $i$-th component $W_i$ of a point $W = (W_k)_{k\in\Nodes}$ in \eqref{eq:Introduction:def-W} as the $i$-th row of a matrix in $\R^{m\times n}$.
Hence points $W \in \WM$ are viewed as row-stochastic matrices with full support, called \textit{assignment matrices}, with \textit{assignment vectors} \eqref{eq:component_i_of_W} as row vectors. 
The barycenter \eqref{eq:def-barycenter-W-product} can then also be expressed as a matrix
\begin{equation}\label{eq:barycenter-WM}
  \BW = \eins_m \BS^\T = \tfrac{1}{n} \eins_m \eins_n^\T.
\end{equation}
The tangent space of $\SM \subset \R^n$ from \eqref{eq:Introduction:def-S} at any point $p \in \SM$ is identified as
\begin{equation}\label{eq:Mechanics_of_Assignment_Flows:tangent-space-S}
  T_p \SM = \{ v \in \R^n \colon \la v, \eins_n\ra = 0\} =: \TS.
\end{equation}
Hence $T_p\SM$ is represented by the same vector subspace $\TS$ of codimension 1, for all $p \in \SM$.
In particular, 
the tangent bundle is trivial
\begin{equation}\label{eq:Mechanics_of_Assignment_Flows:TS}
  T\SM = \SM \times \TS.
\end{equation}
Viewing $\WM$ as an embedded submanifold of $\R^{m\times n}$ by \eqref{eq:Mechanics_of_Assignment_Flows:embedding-W}, we
accordingly identify
\begin{equation}\label{eq:Mechanics_of_Assignment_Flows:tangent-space-W-matrix-embedding}
  T_W \WM = \{ V \in \R^{m\times n} \colon V \eins_n = 0\} =: \TW, \quad \text{for all } W \in \WM \subset \R^{m\times n}.
\end{equation}
With this identification the tangent bundle is also trivial
\begin{equation}\label{eq:Mechanics_of_Assignment_Flows:TW-and-TWdual}
  T\WM = \WM \times \TW.
\end{equation}

\subsection{Assignment Flows} \label{sec:AF} 

Assignment flows are dynamical systems on $\WM$ for inferring probabilistic label assignments that gradually become unambiguous label assignments as $t\to\infty$. These dynamical systems have the form
\begin{equation}\label{eq:Introduction:def-assignment-flow}
  \dot{W}(t) = \ROW_{W(t)}[ F(W(t))], \quad \text{with}\quad W(0) \in \WM,
\end{equation}
where
\begin{equation}\label{eq:affinity-F}
F\colon\mc{W}\to\R^{m\times n} 
\end{equation}
is a smooth function and
\begin{equation}
  \ROW_W \colon \R^{m\times n} \to T_W \WM = \TW , \quad \text{for}\quad W \in \WM,
\end{equation}
is the linear \textit{replicator map} defined componentwise
\begin{subequations}\label{eq:AF-rhs}
\begin{align}
\label{eq:AF-rhs-a} \mc{R}_{W}[F(W)] 
&= \big(\ROS_{W_{i}} F_{i}(W)\big)_{i\in\mc{V}},\quad W\in\mc{W}, \\ 
\intertext{via the \textit{replicator matrix}}
\label{eq:AF-rhs-b}  \ROS_{W_{i}} &= \Diag(W_{i})-W_{i}W_{i}^{\T},\quad 
i\in\mc{V}.
\end{align}
\end{subequations}
The function $F$ couples the evolution of the individual assignment vectors $\dot W_{i},\,i\in\mc{V}$, over the graph, typically by reinforcing tangent directions of similar assignment vectors, and is therefore called \textit{affinity} or \textit{similarity mapping}. Each choice of a similarity mapping $F$ defines a particular assignment flow; see Section~\ref{subsec:admissible-affinity-functions} for a basic instance. Our main result stated in Section \ref{sec:action-functional} characterizes a general \textit{class} of \textit{admissible} similarity mappings $F$.

\subsection{Fisher-Rao Metric}\label{sec:FR-Metric} 

From an information geometric viewpoint \cite{amari2007methods,ay2017information}, the canonical Riemannian structure on $\SM$ is given by the \textit{Fisher-Rao (information) metric}
\begin{equation}\label{eq:def-Fisher-Rao-metric-S}
  g_p(u, v) := \Big\la u, \frac{v}{p} \Big\ra, \quad \text{for all } p \in \SM \;\text{and}\; u, v \in \TS.
\end{equation}
This naturally extends to the product manifold structure of $\WM$ \eqref{eq:Introduction:def-W} via the product metric
\begin{equation}\label{eq:product-Fisher-Rao-metric-W}
  g_W(U, V) := \sum_{i \in [m]} g_{W_i}(U_i, V_i) = \Big\la U, \frac{V}{W}\Big\ra, \quad \text{for all } W \in \WM \;\text{and}\; U, V \in \TW.
\end{equation}
which turns the assignment manifold $\mc{W}$ into a Riemannian manifold.

The orthogonal projection onto $\TS$ and $\TW$, respectively, with respect to the Euclidean inner product are given by 
\begin{subequations}
\begin{align}
  \PTS &\colon \R^n \to \TS, &
  \PTS &:= I_n - \tfrac{1}{n}\eins_n \eins_n \in \R^{n\times n},
  \\
  \PTW &\colon \R^{m\times n} \to \TW, &
  \PTW[A] &:= \big( \PTS A_i\big)_{i\in\Nodes}. \label{eq:def-PTW}
\end{align}
\end{subequations}
Next, we return to the replicator mappings \eqref{eq:AF-rhs}. The linear mapping
\begin{equation}\label{eq:def-repl-matrix}
\ROS_p \colon \R^n \to \TS,\qquad
  \ROS_p = \Diag(p) - pp^\T\in \R^{n\times n}  
\end{equation}
is symmetric
\begin{equation}\label{eq:ROS_symmetric}
  \ROS_p^* = \ROS_p^\T = \ROS_p,
\end{equation}
satisfies the relations
\begin{subequations}\label{eq:Rp-relations}
\begin{align}
\ROS_p &= \ROS_p \PTS = \PTS \ROS_p, 
\label{eq:Rp-relations-a} \\ \label{eq:Rp-relations-b}
\ker(\ROS_p) &= \R\eins_n,
\end{align}
\end{subequations}
and the restriction $\ROS_p|_\TS \colon \TS \to \TS$
to the linear subspace $\TS \subset \R^n$ is a linear isomorphism with inverse given by \cite[Lem.~3.1]{Savarino2019ab}
\begin{equation}\label{eq:inv-repl-matrix}
  (\ROS_p|_{\TS})^{-1} u = \PTS \Diag\big(\tfrac{1}{p}\big) u = \PTS \frac{u}{p}, \quad \text{for all } u \in \TS.
\end{equation}
Likewise, the replicator operator $\ROW_W \colon \R^{m \times n} \to \TW$ satisfies for all $W \in \WM$
\begin{equation}\label{eq:RW-relations}
  \ROW_W = \ROW_W \circ \PTW = \PTW \circ \ROW_W
\end{equation}
and the restriction to the linear subspace $\TW \subset \R^{m\times n}$ is a linear isomorphism with inverse
\begin{equation}\label{eq:inv-RW}
  \big( \ROW_W|_\TW\big)^{-1}[U] = \PTW\Big[\frac{U}{W}\Big], \quad \text{for all } U \in \TW.
\end{equation}
Since all the components $\ROS_{W_i}$ of $\ROW_W$ are symmetric, we have for all $X, Y \in \R^{m\times n}$
\begin{equation}\label{eq:ROW-self-adjoint}
  \la \ROS_{W}[X], Y\ra = \sum_{i \in[m]} \la \ROS_{W_i} X_i, Y_i\ra \overset{\eqref{eq:ROS_symmetric}}{=} \sum_{i \in[m]} \la X_i, \ROS_{W_i} Y_i\ra = \la X, \ROW_W[Y]\ra,
\end{equation}
showing that $\ROW_W$ is self-adjoint $\ROW_W^* = \ROW_W$ with respect to the matrix inner product. There is also a relation between the Fisher-Rao metric and the matrix inner product in terms of the replicator operator. 

\begin{lemma}\label{lem:rel_FRmetric_to_matrix_inner_prod}
  At any point $W \in \WM$, the replicator operator $\ROW_W$ transforms the Riemannian metric into the matrix inner product
  \begin{equation}\label{eq:rel_FRmetric_to_matrix_inner_prod}
    g_W(\ROW_W[U], V) = \la U, V\ra, \quad \text{for all}\quad U, V \in T_W\WM = \TW.
  \end{equation}
\end{lemma}
\begin{proof}
  Using the properties of the replicator operator $\ROW_W$ directly results in
  \begin{subequations}
  \begin{align}
    g_W(\ROW_W[U], V) \overset{\eqref{eq:product-Fisher-Rao-metric-W}}{=}& \Big\la \ROW_W[U], \frac{V}{W} \Big\ra \overset{\eqref{eq:ROW-self-adjoint}}{=} \Big\la U, \ROW_W\Big[\frac{V}{W}\Big]\Big\ra \\
    \overset{\eqref{eq:RW-relations}}{=}& \Big\la U, \ROW_W \circ \PTW\Big[\frac{V}{W}\Big]\Big\ra  
    \overset{\eqref{eq:inv-RW}}{=} \la U, V\ra 
  \end{align}
  \end{subequations}
\end{proof}

\begin{corollary}\label{cor:calculate_ggradJ}
  Let $J \colon \WM \to \R$ be a smooth function and assume there is a smooth map $\Psi \colon \WM \to \R^{m\times n}$ such that the differential of $J$ takes the form
  \begin{equation}
    dJ|_W[V] = \la \Psi(W), V\ra \quad \text{for all } W \in \WM \text{ and } V \in T_W\WM = \TW
  \end{equation}
  with respect to the matrix inner product $\la \cdot, \cdot \ra$. Then, the Riemannian gradient of $J$ is given by
  \begin{equation}
    \ggrad^g J(W) = \ROW_W[\Psi(W)] \quad \text{for all } W \in \WM 
  \end{equation}
\end{corollary}
\begin{proof}
  Let $V \in T_W\WM = \TW$ be arbitrary. As a consequence of Lemma~\ref{lem:rel_FRmetric_to_matrix_inner_prod}, 
  \begin{equation}
    dJ|_W[V] = g_W\big(\ROW_W[\Psi(W)], V\big) \quad \text{for all } V \in T_W\WM = \TW,
  \end{equation}
  with $\ROW_W[\Psi(W)] \in T_W\WM = \TW$. Since this uniquely determines the Riemannian gradient of $J$, the statement follows.
\end{proof}

For functions $J\colon \WM \to \R$ extending onto an open set, the above lemma directly implies a relation between the Riemannian gradient and the usual gradient, a result that is already well known \cite[Prop.~2.2]{ay2017information}. For this, suppose $\wtilde{J} \colon U \to \R$ is a smooth extension of $J$ defined on some open set $U \subset \R^{m\times n}$ containing $\WM$, i.e.~$\wtilde{J}|_\WM = J$. Then, $\Psi(W)$ can be chosen as the usual gradient with respect to the matrix inner product $\partial \wtilde{J}(W) \in \R^{m\times n}$ and the Riemannian gradient of $J$ is given by
\begin{equation}
  \ggrad^g J(W) = \ROW_W[\partial \wtilde{J}(W)], \quad \text{for all } W \in \WM.
\end{equation}

\subsection{The Action Functional}\label{sec:action-functional}

Our main result is summarized in the following theorem. It refers to affinity functions $F \colon \WM \to \R^{m\times n}$ introduced in and discussed after equation \eqref{eq:Introduction:def-assignment-flow}. Applying the identifications $T_W\WM = \TW$ from \eqref{eq:Mechanics_of_Assignment_Flows:tangent-space-W-matrix-embedding} and $T_{F(W)}\R^{m\times n} = \R^{m\times n}$ for every $W \in \WM$ allows to view the differential of $F$ as a linear operator
\begin{equation}
  dF|_W \colon \TW \to \R^{m\times n}.
\end{equation}
The adjoint of $dF|_W$ with respect to the standard matrix inner product \eqref{eq:def-adjoint-lin-op} on $\R^{m\times n}$ and $\TW \subset \R^{m\times n}$ is denoted by
\begin{equation}\label{eq:def-dF-adjoint}
  dF|^*_W \colon \R^{m\times n} \to \TW.
\end{equation}

\begin{theorem}\label{thm:Introduction:main-theorem}
  Let $F\colon \WM \to \R^{m\times n}$ be an affinity map and $W \colon [t_0, t_1] \to \WM$ a solution of the corresponding assignment flow \eqref{eq:Introduction:def-assignment-flow}. Then $W(t)$ is a critical point of the action functional
  \begin{equation}\label{eq:Introduction:var-functional-L-on-W}
    \mc{L}(W) = \int_{t_0}^{t_1} \tfrac{1}{2}\|\dot{W}(t)\|_g^2 + \tfrac{1}{2}\sum_{i \in \Nodes} \vvar_{W_i(t)}\big(F_i(W(t))\big) dt,
  \end{equation}
  if and only if the affinity function $F$ fulfills the condition 
  \begin{equation}\label{eq:Introduction:main-theorem:characterization}
    0 = \ROW_{W(t)} \circ ( dF|_{W(t)} - dF|_{W(t)}^*)\circ \ROW_{W(t)} [F(W(t))]\quad \text{ for }\; t \in [t_0, t_1],
  \end{equation}
  where $dF|_{W(t)}^*$ is the adjoint linear operator of $dF|_{W(t)}$ from \eqref{eq:def-dF-adjoint} and  $\ROW_{W(t)}$ is the replicator operator defined by \eqref{eq:AF-rhs-a}. This condition is equivalent to the Euler-Lagrange equation
  \begin{equation}\label{eq:Introduction:main-theorem:E-L-Eq}
    D_t^g \dot{W}(t) = \tfrac{1}{2} \sum_{i\in\Nodes}\ggradg \vvar_{W_i(t)}\big(F_i(W(t))\big) \quad \text{for }\; t \in [t_0, t_1].
  \end{equation}
\end{theorem}

\begin{remark}

Theorem \ref{thm:Introduction:main-theorem} characterizes the class of affinity functions, in terms of condition \eqref{eq:Introduction:main-theorem:characterization}, for which solutions to the assignment flow equation \eqref{eq:Introduction:def-assignment-flow} are stationary points of the action functional \eqref{eq:Introduction:var-functional-L-on-W} and the Euler-Lagrange equation \eqref{eq:Introduction:main-theorem:E-L-Eq}, respectively. We defer most of the further discussion to Section \ref{sec:Discussion} but mention one important point here. Since every first-order ODE can trivially be described as a special case of the Euler-Lagrange equation of some quadratic Lagrangian it is worth pointing out that the Lagrangian $\mc{L}$ in Theorem \ref{thm:Introduction:main-theorem} is \textit{classical}, that is, of the form \textit{kinetic minus potential energy}. In particular, the potential $-\frac{1}{2}\sum_{i \in \Nodes} \vvar_{W_i(t)}\big(F_i(W(t))\big) $ (note the minus sign) is a non-positive function. Solutions of the assignment flow equation \eqref{eq:Introduction:def-assignment-flow} correspond precisely to those solutions of the Euler-Lagrange equation with energy $0$. Since $0$ is the maximum of the potential this energy value is precisely the Ma\~{n}\'e critical value of this Lagrangian system, see Section \ref{sec:Discussion} for further remarks.
\end{remark}

\vspace{0.2cm}
We proceed with Lemmata to prepare the proof of Theorem \ref{thm:Introduction:main-theorem}.
\begin{lemma}\label{lem:Mechanics_of_Assignment_Flows:relation-sqrd-norm-RO-and-var-on-S}
Let $p \in \SM$ and $f \in \R^{n}$. Then 
\begin{equation}  \label{eq:Var-p-f}
  \|\ROS_p f\|_g^2 = 
  \la f, \ROS_p f\ra = \EE_p[f^2] - (\EE_p[f])^2
  = \vvar_{p}(f).
\end{equation} 
  Thus, for $W \in \WM$ and $F \in \R^{m\times n}$, we have 
\begin{equation}\label{eq:lem-RW-FW-2}
  \|\ROW_W[F]\|_g^2 = \la F, \ROW_W[F]\ra = \sum_{i\in\Nodes} \vvar_{W_i} (F_i).
\end{equation}
\end{lemma}
\begin{proof}
We have
\begin{subequations}
\begin{align}
  \|\ROS_p f\|_g^2 &= g_p(\ROS_p f, \ROS_p f) \overset{\eqref{eq:rel_FRmetric_to_matrix_inner_prod}}{=} \la f, \ROS_p f\ra 
  \overset{\eqref{eq:def-repl-matrix}}{=}
  \la f, p\diamond f - \la p, f\ra p\ra\\
  &= \la f^{\diamond 2}, p\ra - \la f, p\ra^2 = \EE_p[f^2] - (\EE_p[f])^2 = \vvar_p(f).
\end{align}
\end{subequations}
Therefore, it follows
\begin{equation}
  \|\ROW_W[F]\|_g^2 \overset{\eqref{eq:rel_FRmetric_to_matrix_inner_prod}}{=} \la F, \ROW_W[F]\ra \overset{\eqref{eq:AF-rhs-a}}{=} \sum_{i\in[m]} \la F_i, \ROS_{W_i} F_i\ra = \sum_{i\in\Nodes} \vvar_{W_i} (F_i)
\end{equation}
\end{proof}
Next, we compute the differential of the assignment flow vector field \eqref{eq:Introduction:def-assignment-flow} viewed as a mapping
\begin{equation}\label{eq:def-ROW-F}
\ROW[F] \colon \WM \to \TW, \qquad W \mapsto \ROW[F](W):= \ROW_W[F(W)].
\end{equation}
\begin{lemma}\label{lem:Mechanics_of_Assignment_Flows:differential-R[F]}
  With the identifications $T_W\WM = \TW$ and $T_{\ROW_W[F(W)]}\WM = \TW$ due to \eqref{eq:Mechanics_of_Assignment_Flows:tangent-space-W-matrix-embedding}, the differential of the mapping \eqref{eq:def-ROW-F} is a linear map $d\ROW[F]|_W \colon \TW \to \TW$, given by 
  \begin{subequations}
  \begin{align}
    d\ROW[F]|_W[V] 
    &= \ROW_W\circ dF|_W[V] + \mc{B}(W, F(W))[V], \qquad V \in \TW,
    \intertext{where the $i$-th row of the linear map $\mc{B}(W, F(W)) \colon \TW \to \TW$ is defined via matrix multiplication} \label{eq:def-mc-B}
    (\mc{B}(W, F)[V])_i 
    &:= B(W_i, F_i)V_i, \qquad i\in\mc{V},\quad
    W \in \WM,\; F \in \R^{m\times n},
    \\ 
    \intertext{with matrix $B$ given by}
    \label{eq:def-mc-B-component}
    B(p, f) 
    &:= \Diag(f) - \la p, f\ra I_n - p f^\T, \qquad p \in \SM, f \in \R^n.
  \end{align}
  \end{subequations}
\end{lemma}
\begin{proof}
  A short calculation shows $\la B(W_i, F_i(W))V_i, \eins_n\ra = 0$ for all $i \in \Nodes$, that is $\mc{B}(W, X)[V] \in \TW$. Let $\eta \colon (-\veps, \veps) \to \WM$ be a curve with $\eta(0) = W$ and $\dot{\eta}(0) = V$. Keeping in mind $\ROS_p = \Diag(p) - pp^\T$, we obtain for each row vector indexed by $i \in \Nodes$
  \begin{subequations}
  \begin{align}
    \big(d\ROW[F]|_W[V]\big)_i &= \tfrac{d}{dt} \ROS_{\eta_i(t)} F_i(\eta(t))\big|_{t = 0} 
    = \tfrac{d}{dt} \ROS_{\eta_i(t)}\big|_{t=0} F_i(W) +  \ROS_{W_i} \tfrac{d}{dt} F_i(\eta(t))\big|_{t = 0}\\
    &= \big(\Diag(V_i) - V_i W_i^\T - W_i V_i^\T\big) F_i(W) + \big(\ROW_{W} \big[\tfrac{d}{dt} F(\eta(t))\big|_{t = 0}\big]\big)_i\\
    &= \big(\mc{B}(W, F(W))[V]\big)_i + \big(\ROW_W\circ dF|_W[V]\big)_i,
  \end{align}
  \end{subequations}
  where $\Diag(V_i)F_i(W) = \Diag(F_i(W))V_i$ and $V_i^\T F_i(W) = F_i(W)^\T V_i$ was used to obtain the last equality.
\end{proof}

Next, we consider the covariant derivative of a vector fields along a curve $p \colon I \to \SM$, with $I \subset \R$ an interval. Due to $T\SM = \SM \times \TS$, we view a vector field $v(t)$ along $p(t)$ as a map $v \colon I \to \TS$, and also its usual time derivative $\dot{v} \colon I \to \TS$, since $\TS$ is a vector space. 
Defining
\begin{equation}\label{eq:proof-def-Apv}
  A \colon \SM \times \TS \to \TS,\qquad (p, v) \mapsto
  A(p, v):= \frac{v^{\diamond 2}}{p} - \|v\|_g^2 p
\end{equation}
and using the expression from \cite[Eq.~(2.60)]{ay2017information} (with $\alpha$ set to $0$), the covariant derivative $D_t^g v$ of $v$ is related to $\dot{v}$ by
\begin{equation}\label{eq:covariant-derivative-on-S}
  D_t^g v(t) = \dot{v}(t) - \tfrac{1}{2} \frac{(v(t))^{\diamond 2}}{p(t)} + \tfrac{1}{2} \|v(t)\|_g^2 p(t) = \dot{v}(t) - \tfrac{1}{2} A(p(t), v(t)).
\end{equation}
Similarly, as a consequence of $T\WM = \WM \times \TW$, we regard a vector field $V(t)$ along a curve $W \colon I \to \WM$ as a mapping $V \colon I \to \TW$, and likewise $\dot{V} \colon I \to \TW$. Since the covariant derivative on a product manifold equipped with a product metric is the componentwise application of the individual covariant derivatives, the covariant derivative of $V$ on $\WM$ has the form
\begin{equation}\label{eq:Mechanics_of_Assignment_Flows:covariant-derivative-on-W}
  D_t^g V(t) = \dot{V}(t) - \tfrac{1}{2}\mc{A}(W(t), V(t)),
\end{equation}
with $i$-th row of the last term given componentwise by \eqref{eq:proof-def-Apv}
\begin{equation}\label{eq:def-mcA}
  \mc{A} \colon \WM \times \TW \to \TW,\qquad (\mc{A}(W, V))_i = A(W_i, V_i) \quad \text{for all } i\in [m].
\end{equation}
The \textit{acceleration} of a curve $W(t)$ on $\WM$ is the covariant derivative of its velocity vector field $V(t) := \dot{W}(t)$, related to the ordinary time derivative $\dot{V} = \ddot{W}$ in $\R^{m\times n}$ by
\begin{equation}\label{eq:Mechanics_of_Assignment_Flows:acceleration-on-W}
  D_t^g \dot{W}(t) = \ddot{W}(t) - \tfrac{1}{2}\mc{A}(W(t), \dot{W}(t)).
\end{equation}
\begin{lemma}\label{lem:Mechanics_of_Assignment_Flows:cov-der-solution-replicator}
  Suppose $W \colon I \to \WM$ is a solution of the assignment flow \eqref{eq:Introduction:def-assignment-flow}. Then the acceleration of $W(t)$ in terms of the covariant derivative of $\dot{W}(t)$ takes the form
\begin{equation}\label{eq:cov-der-W}
D^g_t \dot{W} = \ROW_{W} \circ dF|_{W} \circ \ROW_{W}[F(W)] + \tfrac{1}{2} \mc{A}\big(W, \ROW_{W}[F(W)]\big).
\end{equation}
\end{lemma}
\begin{proof}
  Since $W(t)$ is a solution of $\dot{W}(t) = \ROW_{W(t)}[F(W(t))]$, the second derivative $\ddot{W}(t) = \frac{d}{dt} \dot{W}(t)$ takes the form (to simplify notation we omit the argument $t$)
  \begin{subequations}
  \begin{align}
    \ddot{W} &= \tfrac{d}{dt}\ROW_{W}[F(W)] = d\ROW[F]|_{W}[\dot{W}] \overset{\text{Lem.~\ref{lem:Mechanics_of_Assignment_Flows:differential-R[F]}}}{=} \ROW_{W}\circ dF|_{W}[\dot{W}] + \mc{B}(W, F(W))[\dot{W}]
    \\
    &= \ROW_{W} \circ dF|_{W}\circ \ROW_{W}[F(W)]
    + \mc{B}(W, F(W))[\ROW_{W}[F(W)]],
  \end{align}
  \end{subequations}
where $\mc{B}$ is defined by \eqref{eq:def-mc-B}. We have $\la f, \ROS_p f\ra = \|\ROS_p f\|_g^2$ by Lemma~\ref{lem:Mechanics_of_Assignment_Flows:relation-sqrd-norm-RO-and-var-on-S} and using \eqref{eq:def-mc-B-component} 
  \begin{subequations}\label{eq:proof-Bpf-Rp}
  \begin{align}
    B(p, f) R_p f &= (f - \la p, f\ra \eins_n) \diamond (\ROS_p f) - \la f, \ROS_p f\ra p\\
    &= \tfrac{1}{p} (\ROS_p f)^{\diamond 2} - \|\ROS_p f\|_g^2 p = A(p, \ROS_p f).
  \end{align}
  \end{subequations}
  This implies $\mc{B}(W, F(W))[\ROW_{W}[F(W)]] = \mc{A}(W, \ROW_W[F(W)])$ and results in the identity
  \begin{equation}
    \ddot{W} = \ROW_{W} \circ dF|_{W} \circ \ROW_{W}[F(W)] + \mc{A}\big(W, \ROW_{W}[F(W)]\big). 
  \end{equation}
  Substituting this expression into \eqref{eq:Mechanics_of_Assignment_Flows:acceleration-on-W} yields \eqref{eq:cov-der-W}.
\end{proof}
As a final preparatory step, we define the potential
\begin{equation}\label{eq:def-G-on-W}
  G \colon \WM \to \R, \qquad G(W) := - \tfrac{1}{2}\|\ROW_W[F(W)]\|_g^2 \overset{\eqref{eq:lem-RW-FW-2}}{=} -\tfrac{1}{2}\sum_{k\in\Nodes} \vvar_{W_k} (F_k(W))
\end{equation}
and compute its Riemannian gradient.

\begin{lemma}\label{lem:Can-Sympl-Struct:Lagrange-Formulation:First-Variation:R-grad-VarF}
  The Riemannian gradient of the potential $G$ from \eqref{eq:def-G-on-W} is given by 
  \begin{equation}\label{eq:grad-GW}
    \ggradg G(W) = -\ROW_W \circ dF|_W^*\circ \ROW_W[F(W)] - \tfrac{1}{2}\mc{A}(W, \ROW_W[F(W)]),\quad\forall W \in \WM,
  \end{equation}
  where $dF|_{W(t)}^*$ is the adjoint linear operator of $dF|_{W(t)}$ from \eqref{eq:def-dF-adjoint}.
\end{lemma}
\begin{proof}
  Let $W \in \WM$. In the following, we derive the expression in \eqref{eq:grad-GW} by applying Corollary~\ref{cor:calculate_ggradJ}. To this end, take any $V \in T_W\WM =\TW$ and let $\eta \colon (-\veps, \veps) \to \WM$ be a curve with $\eta(0) = W$ and $\dot{\eta}(0) = V$. Then 
  \begin{subequations}
  \begin{align}
    dG|_W[V] &= \tfrac{d}{dt}G(\eta(t))\big|_{t=0} \overset{\text{Lem.~\ref{lem:Mechanics_of_Assignment_Flows:relation-sqrd-norm-RO-and-var-on-S}}}{=} -\tfrac{1}{2}\tfrac{d}{dt}\big\la F(\eta(t)), \ROW_{\eta(t)} [F(\eta(t))]\big\ra \big|_{t=0}\\
    &= -\tfrac{1}{2}\big\la \tfrac{d}{dt} F(\eta(t))\big|_{t=0}, \ROW_{\eta(t)} F(\eta(t))\big\ra  -\tfrac{1}{2}\big\la F(\eta(t)), \tfrac{d}{dt}\ROW_{\eta(t)} F(\eta(t))\big|_{t=0}\big\ra\\ 
    &= -\tfrac{1}{2}\big\la dF|_W[V], \ROW_W [F(W)]\big\ra - \tfrac{1}{2}\big\la F(W), d\ROW[F]|_W[V]\big\ra.
  \end{align}
  \end{subequations}
  Using the expression for $d\ROW[F]|_W$ from  Lemma~\ref{lem:Mechanics_of_Assignment_Flows:differential-R[F]} and $\ROW_W^* = \ROW_W$ from \eqref{eq:ROW-self-adjoint}, the second inner product takes the form
  \begin{subequations}
  \begin{align}
    \big\la F(W),& d\ROW[F]|_W[V]\big\ra = \big\la F(W), \ROW_W\circ dF|_W[V]\big\ra + \big\la F(W), \mc{B}(W, F(W))[V]\big\ra\\
    &= \big\la dF|_W^*\circ \ROW_W[F(W)], V\big\ra + \big\la \mc{B}^*(W, F(W))[F(W)], V\big\ra.
  \end{align}
  \end{subequations}
  Substituting back this formula into the above expression for $dG|_W$ together with the expression 
  \begin{equation}
  \big\la dF|_W[V], \ROW_W [F(W)]\big\ra = \big\la V, dF|_W^*\circ \ROW_W [F(W)]\big\ra
  \end{equation}
  for the first inner product, results in
  \begin{subequations}
  \begin{align}
    dG|_W[V] &= \big\la -dF|_W^*\circ \ROW_W[F(W)] - \tfrac{1}{2} \mc{B}^*(W, F(W))[F(W)], V\big\ra \\
    &= \la \Psi(W),V\ra.
  \end{align}
  \end{subequations}
  Due to Corollary~\ref{cor:calculate_ggradJ}, the Riemannian gradient is given by
  \begin{subequations}\label{eq:proof-grad-GW}
  \begin{align}
    \ggrad^g G(W) &= \ROW_W[\Psi(W)]\\
    &= -\ROW_W \circ dF|_W^*\circ \ROW_W[F(W)] - \tfrac{1}{2}\ROW_W[\mc{B}^*(W, F(W))[F(W)]].
  \end{align}
  \end{subequations}
  Regarding the adjoint mapping $\mc{B}^{\ast}$, we have 
  \begin{equation}\label{eq:proof-mcB-ast}
    (\mc{B}^{\ast}(W,F)[U])_i = B(W_i, F_i)^\T U_i \quad \text{for all } i \in \Nodes
  \end{equation}
  and by \eqref{eq:def-mc-B-component}
  \begin{subequations}
  \begin{align}
    B(p,f)R_{p}
    =& (\Diag(f)-\la p,f\ra I_{n}-p f^{\T})(\Diag(p)-p p^{\T})
    \\
    \begin{split}
      =& \Diag(f\diamond p)-\la p,f\ra\Diag(p)-p (f\diamond p)^{\T}\\
      & -(f\diamond p)p^{\T}+ \la p,f\ra p p^{\T} + pp^\T fp^\T
    \end{split}
    \\
    =& (\Diag(p)-p p^{\T})(\Diag(f)-\la p,f\ra I_{n}-f p^{\T})
    \\
    =& R_{p}B(p,f)^{\T}.
  \end{align}
  \end{subequations}
  Thus, by \eqref{eq:proof-Bpf-Rp}, we obtain $R_{p}B(p,f)^{\T} f = A(p,R_{p} f)$ and consequently by the componentwise definitions of $\mc{B}^*$ in \eqref{eq:proof-mcB-ast}, $\mc{A}$ in \eqref{eq:def-mcA} and $\ROW_W$ in \eqref{eq:AF-rhs-a},
  \begin{equation}
  \ROW_W[\mc{B}^*(W, F(W))[F(W)]] = \mc{A}(W, \ROW_W[F(W)]).
  \end{equation}
  Substitution into \eqref{eq:proof-grad-GW} yields \eqref{eq:grad-GW}.
\end{proof}
\begin{proof}[Proof of Theorem \ref{thm:Introduction:main-theorem}] 
Due to Lemma \ref{lem:Mechanics_of_Assignment_Flows:relation-sqrd-norm-RO-and-var-on-S}, the Lagrangian of the action functional 
\eqref{eq:Introduction:var-functional-L-on-W} has the form
\begin{equation}\label{eq:L-by-G}
L(W, V) = \tfrac{1}{2}\|V\|_g^2 - G(W),
\end{equation}
with $G(W)$ defined by \eqref{eq:def-G-on-W}. Therefore, the Euler-Lagrange equation \eqref{eq:Introduction:main-theorem:E-L-Eq} is a direct consequence of Proposition~\ref{prop:Lagrangian-Systems:E-L-Eq-on-Rmflds}. Due to Lemma~\ref{lem:Mechanics_of_Assignment_Flows:cov-der-solution-replicator} and \ref{lem:Can-Sympl-Struct:Lagrange-Formulation:First-Variation:R-grad-VarF}, the expression for the acceleration of $W(t)$ and the Riemannian gradient of $G$ at $W(t)$ both contain the term 
\begin{equation}
  \tfrac{1}{2}\mc{A}(W(t), \ROW_{W(t)}[F(W(t))]) 
\end{equation}
with opposite signs, which yields the relation
\begin{align*}
  D_t^g \dot{W}(t) & + \ggradg G(W(t))
  \\
  &= \ROW_{W(t)}\circ dF|_{W(t)}\circ \ROW_{W(t)} [F(W(t))] - \ROW_{W(t)}\circ dF|_{W(t)}^*\circ \ROW_{W(t)} [F(W(t))]\\
  &= \ROW_{W(t)}\circ (dF|_{W(t)} - dF|_{W(t)}^*)\circ \ROW_{W(t)}[F(W(t))].
\end{align*}
As a consequence, the characterization of $F$ in \eqref{eq:Introduction:main-theorem:characterization} is equivalent to the Euler-Lagrange equation \eqref{eq:Introduction:main-theorem:E-L-Eq} and by Proposition~\ref{prop:Lagrangian-Systems:E-L-Eq-on-Rmflds} equivalent to $W(t)$ being a critical point of the action functional.
\end{proof}
%

\section{Discussion}
\label{sec:Discussion}

\subsection{Some Implications of Theorem \ref{thm:Introduction:main-theorem}}
We discuss in the this section various properties and  consequences of Theorem \ref{thm:Introduction:main-theorem}.

\subsubsection{Ma\~n\'e critical value}
In his influential work \cite{mane1997} Ma\~n\'e introduced \textit{critical values} which should be interpreted as energy levels that mark important dynamical and geometric changes for the Euler–Lagrange flow, see \cite{abbo2013} for a nice introduction. Dynamical properties at energies being equal to a Ma\~n\'e critical value are often times hard to analyze. In general, there are various related Ma\~n\'e critical values, however for classical Lagrangians such as $L$, e.g.~\eqref{eq:L-by-G}, in Theorem \ref{thm:Introduction:main-theorem} all of them agree and equal the maximum of the potential. As pointed out before the potential part of the Lagrangian $L$ is $G(W) = -\frac{1}{2}\sum_{i \in \Nodes} \vvar_{W_i(t)}\big(F_i(W(t))\big)$ which has $0$ as maximum. At the same time solutions to the assignment flow equation \eqref{eq:Introduction:def-assignment-flow} are precisely the solutions to the Euler-Lagrange equation \eqref{eq:Introduction:main-theorem:E-L-Eq} of energy $0$, i.e.~at the Ma\~n\'e critical value of $\mc{L}$. 

In the following, basic properties of the set of Ma\~n\'e critical points on $\WM$
\begin{equation}\label{eq:def-Mane-crit-points}
  \Mcrit := \argmax_{W \in \WM} G(W) = G^{-1}(0)
\end{equation}
are investigated and summarized in Proposition~\ref{prop:Mcrit-measure-zero-Mreg-dense-submfld}.
Subsequently, based on a result from geometric mechanics, Proposition~\ref{prop:a-flow-is-geodesic} shows that integral curves of the assignment flow that are critical points of the action functional $\mc{L}$ in Theorem~\ref{thm:Introduction:main-theorem} and start in the complement
\begin{equation}\label{eq:def-complement-of-Mane-circ-points}
  \Mreg := \WM \setminus \Mcrit
\end{equation}
are actually be reparametrized geodesics of the so called \textit{Jacobi metric} introduced below.\\

By Lemma~\ref{lem:Mechanics_of_Assignment_Flows:relation-sqrd-norm-RO-and-var-on-S}, we have 
\begin{equation}\label{eq:Mcrit-points-are-AF-equilibrium-points}
  0 = G(W) \overset{\eqref{eq:def-G-on-W}}{=} -\tfrac{1}{2}\|\ROW_W[F(W)]\|_g^2 \quad \Leftrightarrow \quad 0 = \ROW_W[F(W)],
\end{equation}
that is the potential assumes its maximum at $W$ if and only if $W$ is an equilibrium point of the assignment flow \eqref{eq:Introduction:def-assignment-flow}. 
Due to $\ROW_W|_\TW$ being a linear isomorphism by \eqref{eq:inv-RW}, we further obtain
\begin{equation}\label{eq:Mcrit-points-alt-characterization}
  0 = \ROW_W[F(W)] \overset{\eqref{eq:RW-relations}}{=} \ROW_W|_\TW \circ \PTW[F(W)] \quad \Leftrightarrow \quad 0 = \PTW[F(W)].
\end{equation}
Thus, we need to consider the zero set of the smooth map 
\begin{equation}
  \PTW\circ F \colon \WM \to \TW.
\end{equation}
We restrict our analysis to affinity functions $F$ for which the differential $d(\PTW \circ F)$ has constant rank on $\WM$, in the following denoted by $r$. To avoid the trivial case $\PTW \circ F \equiv \mrm{const}$ we further restrict to the case $r \geq 1$. A basic instance of this case is given in Section~\ref{subsec:admissible-affinity-functions} with $F$ being a linear map.

Due to the Constant-Rank Level Set Theorem~\cite[Thm.~5.12]{lee2013smooth}, the zero set 
\begin{equation}\label{eq:Mane-crit-points-characterizations}
  (\PTW\circ F)^{-1}(0) = G^{-1}(0) = \Mcrit \subset \WM
\end{equation}
is a properly embedded submanifold of $\WM$ with dimension 
\begin{equation}\label{eq:dim-of-Mane-crit-points}
  \dim(\Mcrit) = \dim(\WM) - r \leq \dim(\WM) - 1.
\end{equation}
Since the dimension of $\Mcrit$ is strictly less than $\dim(\WM)$, it is a submanifold with measure zero in $\WM$  \cite[Cor.~6.12]{lee2013smooth}. Therefore, the complement $\Mreg$ \eqref{eq:def-complement-of-Mane-circ-points}, that is the set of points $W$ with $G(W) < 0$, is a dense (\cite[Prop.~6.8]{lee2013smooth}) subset of $\WM$. According to \cite[Prop.~5.5]{lee2013smooth}, being properly embedded in $\WM$ is equivalent to being a closed subset of $\WM$ (in the subspace topology). Thus, $\Mreg$ is an open subset of $\WM$ and consequently also a submanifold. Overall we have proven the following statement.

\begin{proposition}\label{prop:Mcrit-measure-zero-Mreg-dense-submfld}
  If the differential $d(\PTW \circ F)$ has constant rank $r\geq 1$ on $\WM$, then the set $\Mcrit$ of  Ma\~n\'e critical points \eqref{eq:def-Mane-crit-points} is a submanifold of $\WM$ with measure zero and its complement $\Mreg \subset \WM$ \eqref{eq:def-complement-of-Mane-circ-points} is an open and dense subset.
\end{proposition}

Equipped with this result, we are now able to characterize solutions of the assignment flow \eqref{eq:Introduction:def-assignment-flow} starting in $\Mreg$ as reparametrized geodesics.

\begin{definition}
  \textbf{(\cite[Def.~3.7.6]{abraham1987foundations}).}
  Let $h$ be a Riemannian metric on $M$ and $G \colon M \to \R$ a potential. Assume $C$ is a constant such that $G(x) < C$ holds for all $x \in M$. Then the \textit{Jacobi metric} is defined by
  \begin{equation}
    h_C:= (C - G)h.
  \end{equation}
\end{definition}

\begin{theorem}\label{thm:geodesics-of-Jacobi-metric}
  \textbf{(\cite[Thm.~3.7.7]{abraham1987foundations}).}
  Up to reparametrization, the base integral curves of the Lagrangian $L(x, v) = \frac{1}{2}\|v\|_h^2 - G(x)$ with energy $E_0$ are the same as geodesics of the Jacobi metric $h_{E_0}$ with energy $1$.
\end{theorem}

Since $G < 0$ on $\Mreg$, we restrict our investigation to the Riemannian submanifold $(\Mreg, g|_\Mreg)$ and set $C:= 0$, resulting in the Jacobi metric $h_0 = (-G) g|_\Mreg$ of the form 
\begin{equation}\label{eq:energy-0-Jacobi-metric}
  (h_0)_W \overset{\eqref{eq:def-G-on-W}}{=} \tfrac{1}{2}\sum_{k\in\Nodes} \vvar_{W_k} (F_k(W))\ g_W \quad \text{for any point } W \in \Mreg.
\end{equation}
Now, let $W(t)$ be an integral curve of the assignment flow \eqref{eq:Introduction:def-assignment-flow}. If the initial value $W(0)$ lies in $\Mreg$, then the entire integral curve $W(t)$ remains in $\Mreg$. This is a consequence of Ma\~n\'e critical points being equilibrium points by \eqref{eq:Mcrit-points-are-AF-equilibrium-points} and the fact that the assignment flow is a first-order ODE. If additionally $W(t)$ is a critical point of the action functional $\mc{L}$ from Theorem~\ref{thm:Introduction:main-theorem}, then $W(t)$ is a base integral curve with energy $E_0 = 0$. Thus,  Theorem~\ref{thm:geodesics-of-Jacobi-metric} directly implies the following statement.

\begin{proposition}\label{prop:a-flow-is-geodesic}
  Let $W(t)$ be an integral curve of the assignment flow \eqref{eq:Introduction:def-assignment-flow}. If $W(t)$ is a critical point of the action function $\mc{L}$ in Theorem~\ref{thm:Introduction:main-theorem} with initial value $W(0)\in \Mreg$, then, up to reparametrization, $W(t)$ is a geodesic of the Jacobi metric \eqref{eq:energy-0-Jacobi-metric}.
\end{proposition}
\begin{remark}
  It is important to note that the previous statement is only true for solutions of the assignment flow, which is a first-order ODE. A general solution of the second-order ODE Euler-Lagrange equation \eqref{eq:Introduction:main-theorem:E-L-Eq} might leave $\Mreg$ in finite time and cross the set $\Mcrit$.
\end{remark}

In the next section, we directly determine the set $\Mcrit$ for the a basic instance of an assignment flow.

\subsubsection{Admissible Affinity Functions}\label{subsec:admissible-affinity-functions}

Condition \eqref{eq:Introduction:main-theorem:characterization} characterizes affinity functions $F$ for Theorem \ref{thm:Introduction:main-theorem} to hold. We contrast this condition with a simple affinity function used in prior work and directly determine the corresponding set $\Mcrit$ of Ma\~n\'e critical points.

The recent paper \cite[Proposition 3.6]{Savarino2019ab} introduced a reparametrization, called $S$-\textit{flow}, of the original assignment flow formulation of \cite{Astroem2017}. The distance information between each data point $f_i \in \mc{F}$ and the labels $f^*_j \in \mc{F}$ is collected in the \textit{data matrix}
\begin{equation}
  D \in \R^{m\times n}\quad \text{with}\quad D_{ij} = d_\mc{F}(f_i, f^*_j), \quad \text{for } i \in [m], j \in [n],
\end{equation}
where $d_\mc{F}$ is the metric introduced in Section~\ref{sec:AssManifold}. Intuitively it represents how well each data point is represented by the labels. For a nonnegative averaging matrix 
\begin{equation}\label{eq:Om-properties}
  \Om \in \R^{m\times m} \quad \text{with}\quad \Om \geq 0 \quad \text{and}\quad \Om \eins_m = \eins_m,
\end{equation}
the $S$-flow equations read
\begin{subequations}\label{eq:S-flow}
\begin{align}
  \dot S &= \ROW_{S}[\Omega S], &
  S(0) &= \exp_{\eins_{\mc{W}}}(-\Omega D),
  \label{eq:S-flow-S} \\ \label{eq:S-flow-W}
  \dot W &= \ROW_{W}[S], &
  W(0) &= \eins_{\mc{W}},
\end{align}
\end{subequations}
where the so-called \textit{lifting map}
\begin{subequations}
\begin{align}
  \exp_{W}\colon\TW\to\mc{W},\qquad
  \exp_{W} &= \Exp_{W}\circ \ROW_{W}, \quad W\in\mc{W},
  \\
  \big(\exp_{W}(V)\big)_{i}
  &= \frac{W_{i}\diamond e^{V_{i}}}{\la W_{i},e^{V_{i}}\ra},\quad i\in[m],\quad W\in\mc{W},\; V\in\TW
\end{align}
\end{subequations}
is the composition of the mapping \eqref{eq:AF-rhs-a} and the exponential map $\Exp$ of $(\mc{W},g)$ with respect to the so-called e-connection of information geometry \cite{amari2007methods}. Note that both solutions $S(t), W(t)$ evolve on $\mc{W}$ and that $W(t)$ depends on $S(t)$ but not vice versa. Hence we focus on the system \eqref{eq:S-flow-S} and the specific affinity function given by matrix multiplication
\begin{equation}\label{eq:F-S-flow}
  F(S) = \Omega S.
\end{equation}
The differential of $F$ at $S \in \WM$ is therefore also given by matrix multiplication
\begin{equation}\label{eq:dF-S-flow}
  dF|_{S}[V]=\Omega V,
\end{equation}
that is condition \eqref{eq:Introduction:main-theorem:characterization} holds in particular if $\Omega=\Omega^{\T}$ is symmetric. This assumption was adopted in \cite{Savarino2019ab} and in a slightly more general form also in \cite{Zern:2020aa}.

Next, we determine the set $\Mcrit$ of Ma\~n\'e critical points \eqref{eq:def-Mane-crit-points} based on the condition on the right-hand side of \eqref{eq:Mcrit-points-alt-characterization}. A basic calculation using the properties of $\PTW$ and $\Om$ shows that these two linear operators commute, resulting in
\begin{equation}\label{eq:PTWoF_expressed_with_Om}
  \PTW[F(W)] = \PTW[\Om W] = \Om\PTW[W] = \Om \big( W - \BW\big) \quad \text{for all } W \in \WM.
\end{equation}
Since the corresponding differential is just matrix multiplication independent of $W \in \WM$
\begin{equation}\label{eq:dPTWoF_expressed_with_Om}
  d(\PTW \circ F)|_W[V] = \PTW[\Om V] = \Om V \quad \text{for all } V \in \TW,
\end{equation}
the rank $r$ of $d(\PTW \circ F)$ is constant. For Proposition~\ref{prop:Mcrit-measure-zero-Mreg-dense-submfld} to hold, we need to check that the rank satisfies $r\geq 1$. For this, denote the corresponding kernel of  \eqref{eq:dPTWoF_expressed_with_Om} by
\begin{equation}
  \Sigma_\Om := \ker\big(d(\PTW \circ F)\big) = \{ V \in \TW\ |\ \Om V = 0\}.
\end{equation}

\begin{lemma}\label{lem:dim-Sigma-Om}
  $\dim(\Sigma_\Om) = (n-1)\dim(\ker(\Om))$ and therefore the rank of $d(\PTW \circ F)$ on $\WM$ is $r = (n-1) \rank(\Om)$.
\end{lemma}
\begin{proof}
  Denote the standard basis of $\R^n$ by $e_1, \ldots, e_n$. A basis for $\TS$ \eqref{eq:Mechanics_of_Assignment_Flows:tangent-space-S} is then given by 
  \begin{equation}
    b_i := e_i - e_n, \quad \text{for } i \in [n-1].
  \end{equation}
  Furthermore, set $K:= \dim(\ker(\Om))$ and let $a_1, \ldots, a_K$ be a basis of $\ker(\Om) \subset \R^m$. Then, for every $k \in [K]$ and $i \in [n-1]$
  \begin{equation}
    a_k b_i^\T \eins_n = a_k \la b_i, \eins_n\ra = 0\quad \text{and}\quad \Om a_kb_i^\T = 0,
  \end{equation}
  showing that $a_k b_i^\T \in \Sigma_\Om$.
  As all the $a_k$ and $b_i$ are each linear independent, so are their outer products $a_k b_i^\T$ for all $k \in[K]$ and $i \in [n-1]$. 
  Now, let $V \in \Sigma_\Om$ be arbitrary. Writing $V$ as $V = \sum_{i \in [n]} V e_i e_i^\T$ we obtain
  \begin{equation}
    0 \overset{\eqref{eq:Mechanics_of_Assignment_Flows:tangent-space-W-matrix-embedding}}{=} V\eins_n = \sum_{i\in[n]} Ve_i \quad \Leftrightarrow \quad Ve_n = -\sum_{i\in [n-1]} Ve_i,
  \end{equation}
  which in turn shows that $V$ can be expressed in terms of the basis $b_i$ as $V = \sum_{i\in[n-1]} V e_i b^\T_i$.
  On the other hand, the $i$-th column of $V$, given by $Ve_i$, fulfills $\Om Ve_i = 0$ and can be expressed as
  $Ve_i = \sum_{k\in[K]} \lambda_{ki} a_k$, with coefficients $\lambda_{ki} \in \R$.
  Putting everything together results in
  $V = \sum_{i\in[n-1]} V e_i b^\T_i = \sum_{i\in[n-1]}\sum_{k\in[K]} \lambda_{ki} a_kb^\T_i$, 
  showing that all the $a_k b_i^\T$ are indeed a basis for $\Sigma_\Om$. As a result, the formulas for $\dim(\Sigma_\Om)$ and the rank
  \begin{equation}
    r = \dim(\TW) - \dim(\Sigma_\Om) = (n-1)m - (n-1)\dim(\ker(\Om)) = (n-1)\rank(\Om)
  \end{equation}
  follow.
\end{proof}
As a consequence of $\rank(\Om) \geq 1$  by \eqref{eq:Om-properties}, a lower bound on the rank $r$ of $d(\PTW\circ F)$ is given by $r \geq n-1 \geq 1$. Therefore, Proposition~\ref{prop:Mcrit-measure-zero-Mreg-dense-submfld} applies and $\Mcrit$ for the $S$-flow is a submanifold of $\WM$ with measure zero. The expression of $\PTW\circ F$ in terms of $\Om$ from \eqref{eq:PTWoF_expressed_with_Om} and the fact that $W -\BW$ lies in $\TW$ for all $W \in \WM$ allow to explicit characterization $\Mcrit$ as an affine subspace
\begin{equation}
   \Mcrit \overset{\eqref{eq:Mane-crit-points-characterizations}}{=} (\PTW \circ F)^{-1}(0) \overset{\eqref{eq:PTWoF_expressed_with_Om}}{=} (\BW + \Sigma_\Om) \cap \WM.
\end{equation}
with dimension
\begin{equation}
  \dim(\Mcrit) \overset{\eqref{eq:dim-of-Mane-crit-points}}{=} \dim(\WM) - r \overset{\text{Lem.~\ref{lem:dim-Sigma-Om}}}{=} \dim(\Sigma_\Om) \leq (m-1)(n-1).
\end{equation}
As $\Om$ is assumed to be given, $\Mcrit$ can explicitly be constructed after a basis for $\ker(\Om)$ has been calculated. Therefore, we are able to check if $S(0) \notin \Mcrit$, in which case the corresponding integral curve $S(t)$ of the $S$-flow \eqref{eq:S-flow-S} would be a reparametrized geodesic for the Jacobi metric \eqref{eq:energy-0-Jacobi-metric} with energy $E_0 = 0$, according to Theorem~\ref{thm:geodesics-of-Jacobi-metric}.

We conclude this section with another observation that should stimulate future work. 
Under the afore-mentioned symmetry assumption, a continuous-domain approach was studied in \cite{Savarino2019ab} corresponding to \eqref{eq:S-flow} at `spatial scale zero'. The latter means to consider only parameter matrices $\Omega$ in \eqref{eq:S-flow-S} whose sparse row vectors $\Omega_{i}$ encode \textit{nearest-neighbor} interactions of $S_{i}$ and $\{S_{k}\colon k\sim i\}$ on an underlying regular grid graph, and to consider the right-hand side of \eqref{eq:S-flow-S} as discretized Riemannian gradient of a continuous-domain variational approach with \textit{pointwise} defined variables. Specifically, replacing $i\in\mc{V}$ by locations $x\in U\subset\R^{d}$, the vector field $S \colon \mc{V} \to \SM$, $i \mapsto S_i$, becomes a simplex-valued vector field $S \colon U \to \SM$, $x \mapsto S(x)$, that has to solve a variational inequality. Besides analyzing existence of a minimizer in a suitable function space and a corresponding dedicated numerical algorithm, a heuristically (under too strong regularity assumptions) derived partial differential equation was presented that is supposed to characterize any minimizer $S^{\ast}$ and reads
\begin{equation}\label{eq:S-PDE}
  R_{S^{\ast}}(-\Delta S^{\ast}-\alpha S^{\ast}) = 0,
\end{equation}
where $R_{S^{\ast}}$ applies pointwise $R_{S^{\ast}(x)}$ to the vector $(-\Delta S^{\ast}-\alpha S^{\ast})(x)$ at every $x\in\Omega$, in the same way as the mapping $\ROW_{W}$ defined by \eqref{eq:AF-rhs-a} amounts to applying the mappings \eqref{eq:AF-rhs-b} at every vertex $i\in\mc{V}$.

From this viewpoint, condition \eqref{eq:Introduction:main-theorem:characterization},
\begin{equation}
  0 = \ROW_{W(t)}\circ (dF|_{W(t)}-dF|^{\ast}_{W(t)})\circ\ROW_{W(t)}[F(W(t))],
\end{equation}
that was shown to be equivalent to the Euler-Lagrange equation \eqref{eq:Introduction:main-theorem:E-L-Eq}, should become the spatially-\textit{discrete} but \textit{nonlocal} analogon of \eqref{eq:S-PDE} in the limit $t\to\infty$. We leave the exploration of this observation for future work.

\subsubsection{Geometric Integration Versus Optimization}\label{sec:integration-vs-optimization}

In contrast to classical approaches of the labeling problem, the presented dynamical geometric formulation does not merely rely on finding maximizers of a task specific objective function, but instead solely depends on the Lagrangian dynamics governing the inference process. In the following, this is discussed in more detail.

Classical formulations of image labeling \cite{Kappes:2015aa} are usually formulated as minimization problems of (preferably convex) functions $\min_X J(X)$, where global minimizers are associated with meaningful label assignments. As a consequence, the minimizers themselves are the solution of the labeling problem, independent of any specific optimization strategy used to find or approximate them.

In \cite[Pro.~3.9, Prop.~3.10]{Savarino2019ab} it was shown that if the weight matrix $\Om$ is symmetric $\Om=\Om^{\T}$, then the above mentioned $S$-\textit{flow} \eqref{eq:S-flow-S} is actually a Riemannian gradient ascent flow with respect to the function
\begin{equation}\label{eq:integration-vs-optimization:S-flow-potential-J-def}
  J(S) = \tfrac{1}{2}\la S, \Om S\ra = \tfrac{1}{2}\|S\|_2^2 - \tfrac{1}{4}\sum_{i\in\Nodes}\sum_{j\in\nhood_i} \Om_{ij} \|S_i - S_j\|_2^2.
\end{equation}
Similar to the continuous case in \cite[Prop.~4.2]{Savarino2019ab}, it can be shown that the global maximizers of $J$ are spatially constant assignments, i.e.~every node in the graph has the same label. This can directly be seen from the right-hand side expression for $J$ in \eqref{eq:integration-vs-optimization:S-flow-potential-J-def}. In order for $J$ to obtain its supremum, the first term $\| S\|_2^2$ needs to be maximal, which happens precisely if every $S_i$ is one of the standard basis vectors, and the second term $\sum_{i\in\Nodes}\sum_{j\in\nhood_i} \Om_{ij} \|S_i - S_j\|_2^2$ needs to be minimal (zero), which happens precisely if all the $S_i$ have the same value at all nodes $i \in \Nodes$, that is $S$ is spatially constant.

Therefore, in contrast to the above mentioned classical methods, we are \textit{not} interested in maximizers of the function $J$, as they generally do not represent meaningful assignments. Indeed, any nontrivial assignment the $S$-flow $S(t)$ converges to (which experimentally happens \cite{Zern:2020aa, Savarino2019ab}) \textit{cannot} be a maximizer of $J$. Rather, the integral curves themselves, that is the inference process governed by the spatially coupled replicator dynamics, is the crucial element responsible for producing meaningful label assignments as limit points. This highlights the importance of the Lagrangian mechanical viewpoint of the assignment.

\subsection{Directly Related Work}\label{sec:RK18}

In \cite[Thm.~2.1]{raju2018variational}, the authors claim that all \textit{uncoupled} equations of the form $\dot{p} = \ROS_p F(p)$, on a single simplex $p(t) \in \SM$, satisfy the Euler-Lagrange equation associated with the cost functional 
\begin{equation}\label{eq:Introduction:def-functional-L-on-S}
  \mc{L}(p) := \int_{t_0}^{t_1} \tfrac{1}{2}\|\dot{p}(t)\|_g^2 + \tfrac{1}{2}\|\ROS_{p(t)} F(p(t))\|_g^2 dt
  \quad \text{for curves } p \colon [t_0, t_1] \to \SM.
\end{equation}
In our present paper, we derive a more general result (Theorem \ref{thm:Introduction:main-theorem}) for a system \eqref{eq:AF} of \textit{coupled} equations from the viewpoint of geometric mechanics on manifolds, of  which \eqref{eq:Introduction:def-functional-L-on-S} is a (very) special case. In particular, we derive a necessary condition \eqref{eq:Introduction:main-theorem:characterization} that is missing in \cite{raju2018variational}, which any affinity function $F$ has to satisfy for the assertion of Theorem \ref{thm:Introduction:main-theorem} to hold. This latter result yields an interpretation of stationary points of the action function as solutions of the Euler-Lagrange equation \eqref{eq:Introduction:main-theorem:E-L-Eq}.

It can be shown that in the case of $n = 2$ labels, any fitness function $F$ indeed fulfills condition \eqref{eq:Introduction:main-theorem:characterization} and therefore also the Euler-Lagrange equation.
However, for $n > 2$ labels this is no longer true, as the following counterexample demonstrates.

Suppose we have more than two labels, i.e.~$n>2$, and first consider the case of $m = |\Nodes| = 1$ nodes, that is an \textit{uncoupled} replicator equation on a \textit{single} simplex. Define the matrix $F := e_2 e_1^\T$, where $e_i$ are the standard basis vectors of $\R^n$. Thus, the affinity function is a linear map 
\begin{equation}\label{eq:counterexample-def-F}
  F\colon \SM \to \R^n, \quad p = (p^1, \ldots, p^n)^\T \mapsto Fp = p^1 e_2
\end{equation}
fulfilling $d F_p = F$ and $d F_p^* = F^\T$. A short calculation using the relation $\ROS_p e_i = p^i (e_i - p)$ (Einstein summation convention is \textit{not} used) shows that the first coordinate of condition \eqref{eq:Introduction:main-theorem:characterization} takes the form
\begin{equation*}
  \big(\ROS_p (F - F^\T)\ROS_p F p\big)^1 = - (p^1)^2 p^2 (1 - p^1 - p^2) \neq 0,\quad \text{for all } p \in \SM.
\end{equation*}
This example generalizes to the case $m > 1$ by defining the linear affinity function $\mc{F}[W]$ componentwise by $(\mc{F}[W])_i := F W_i,\,i\in[m]$.

\subsection{Lagrangian and Hamiltonian Point of View}
Theorem \ref{thm:Introduction:main-theorem} rests upon the representation of the assignment flow as a Lagrangian mechanical system of the form kinetic minus potential energy \eqref{eq:L-by-G}, as summarized in Section \ref{sec:Mechanics-Riemannian}. Due to this specific form, Proposition~\ref{prop:Lagrangian-Systems:E-L-Eq-on-Rmflds} can be applied to characterize critical points of the action functional $\mc{L}$ from Theorem~\ref{thm:Introduction:main-theorem} as solutions to the Euler-Lagrange equation \eqref{eq:Introduction:main-theorem:E-L-Eq}, which in turn allows to derive condition \eqref{eq:Introduction:main-theorem:characterization}. 

For general Lagrangians, however, Proposition~\ref{prop:Lagrangian-Systems:E-L-Eq-on-Rmflds} is not applicable and critical points of the action functional are characterized as integral curves of the Lagrangian vector field $X_E$ as detailed in Section~\ref{sec:Lagrangian:Lagrangian-Systems}. Since Lagrangians of the form kinetic minus potential energy \eqref{eq:Lagrangian-Systems:L-Riemann-mfld} are hyperregular, the representation as Hamiltonian system via the Legendre transformation $\F L$ is an equivalent alternative. As mentioned in Section~\ref{sec:Legendre-transform}, the energy $E\colon T\WM \to \R$, the Hamiltonian $H \colon T^*\WM \to \R$ and their corresponding vector fields $X_E$ on $T\WM$ and $X_H$ on $T^*\WM$ are related via
\begin{equation}
  E = H \circ \F L \quad \text{and}\quad X_E = (\F L)_*^{-1} X_H.
\end{equation}
To obtain interpretable explicit formulas, it will be more convenient to work on $T\WM$ instead of $T^*\WM$. In the following, we derive an explicit expression for the Lagrangian vector field $X_E$ and relate its corresponding integral curves to the Euler-Lagrange equation \eqref{eq:Introduction:main-theorem:E-L-Eq} of Theorem~\ref{thm:Introduction:main-theorem}. Because $X_E$ is the symplectic gradient of the energy $E$ with respect to the Lagrangian form $\om_L$, see \eqref{eq:Lagrangian-two-form-general-def}, we first calculate an alternative formula for $\om_L$ in terms of the Fisher-Rao metric. For this we exploit the fact that the assignment manifold is a so called \textit{Hessian manifold} \cite{shima1997geometry}, that is in suitable coordinates the Fisher-Rao metric is the Hessian of a convex function.

Since $T\WM = \WM \times \TW$ \eqref{eq:Mechanics_of_Assignment_Flows:TW-and-TWdual} is trivial, the tangent space of $T\WM$ at any point $(W, V) \in T\WM$ can be identified with the vector space
\begin{equation}
  T_{(W, V)} T\WM = \TW \times \TW.
\end{equation}
With this identification, the Lagrangian two-form $\om_L$ has the following simple expression.

\begin{lemma}\label{lem:Lagriangian-two-form-on-TW}
  Let $(W, V) \in T\WM$ and $A = (A', A''), B = (B', B'') \in T_{(W, V)}T\WM = \TW \times \TW$. Then the Lagrangian two-form can be expressed via the Fisher-Rao metric as
  \begin{equation}
    \om_L|_{(W, V)} \big(A, B\big) = g_W(A', B'') - g_W(A'', B').
  \end{equation}
\end{lemma}
\begin{proof}
  In the following, if $\vphi$ is a real valued function on $\SM$ or $\WM$, then its coordinate representation is denoted by $\what{\vphi}$. A global chart on $\SM$ is given by $\eta_\SM \colon \SM \to \R^{n-1}$ with $p \mapsto \eta_\SM(p) = (p^1, \ldots, p^{n-1})$. It is a standard result from information geometry \cite{amari2010information} that the
  \textit{negative entropy} $\vphi$, a smooth convex function on $\SM$ defined by
  \begin{equation}
    \vphi \colon \SM \to \R, \quad p \mapsto \sum_{i \in [n]} p^i \log(p^i) = \la p, \log(p)\ra,
  \end{equation}
  induces the Fisher-Rao metric in coordinates $\eta_\SM$, denoted by $(g^\SM_{ij})$, as the Hessian of $\what{\vphi}$
  \begin{equation}
    g^\SM_{ij}(p) = \frac{\partial^2 \what{\vphi}}{\partial p^i \partial p^j} (p^1, \ldots, p^{n-1}) 
    \quad \text{for all } i, j \in [n-1]. 
  \end{equation}
  Thus, a single simplex $\SM$ has the structure of a Hessian manifold \cite{shima1997geometry}. As a global chart of the product manifold $\WM = \prod_{i\in[m]} \SM$ we take the product chart $\eta_\WM \colon \WM \to \R^{m(n-1)}$ with $W \mapsto \eta_\WM (W) = (\eta_\SM(W_1), \ldots, \eta_\SM(W_m)) =  (x^1, \ldots, x^{m(n-1)}) = x$, where each $W_i$ lies in $\SM$ for all $i \in [m]$. 
  Define the \textit{accumulated negative entropy} by
  \begin{equation}
    \vphi_\mrm{acc} \colon \WM \to \R, \quad W \mapsto \sum_{i\in[m]} \vphi(W_i).
  \end{equation}
  and let $(g^\WM_{ij})$ denote the representation of the product Fisher-Rao metric \eqref{eq:product-Fisher-Rao-metric-W} on $\WM$ in coordinates $\eta_\WM$. Since $\vphi_\mrm{acc}$ separates over the product structure of $\WM$, the accumulated negative entropy also induces the product Riemannian metric in the chart $\eta_\WM$ 
  \begin{equation}\label{eq:Lagriangian-two-form-on-TW:g-as-second-derivative}
    g^\WM_{ij}(x) = \frac{\partial^2 \what{\vphi}_\mrm{acc}}{\partial x^i \partial x^j} (x),
  \end{equation}
  equipping also the assignment manifold with the structure of a Hessian manifold \cite{shima1997geometry}.

  Now, take an arbitrary point $(W, V) \in T\WM = \WM \times \TW$ and let $(x, v)$ be the corresponding coordinates with respect to the chart $\eta_\WM$. According to \cite[Prop.~3.5.6]{abraham1987foundations}, the Lagrangian two-form $\om_L$ \eqref{eq:Lagrangian-two-form-general-def} in coordinates is given by
  \begin{equation}\label{eq:Lagriangian-two-form-on-TW:local-expr-general-theta_L}
    \om_L = \sum_{i, j} \bigg( \frac{\partial^2 L}{\partial v^i \partial x^j } dx^i \wedge dx^j 
      + \frac{\partial^2 L}{\partial v^i \partial v^j} dx^i\wedge dv^j\bigg).
  \end{equation}
  Since the coordinate expression of the Lagrangian~\eqref{eq:L-by-G} is $L(x, v) = \frac{1}{2}\sum_{i, j} g^\WM_{ij}v^i v^j - G(x)$, the second-order derivatives are 
  \begin{equation}
    \frac{\partial^2 L}{\partial v^i \partial x^j} = \sum_k\frac{\partial g^\WM_{ik}}{\partial x^j}v^k
    \quad \text{and}\quad
    \frac{\partial^2 L}{\partial v^i \partial v^j} = g^\WM_{ij}.
  \end{equation}
  Plugging these expressions into \eqref{eq:Lagriangian-two-form-on-TW:local-expr-general-theta_L} and rearranging the first sum using $dx^j\wedge dx^i = - dx^i \wedge dx^j$ yields
  \begin{equation}\label{eq:Lagriangian-two-form-on-TW:local-expr-Riemannian-theta_L}
    \om_L = \sum_{i < j} \sum_{k}\Big(\frac{\partial g^\WM_{ik}}{\partial x^j} - \frac{\partial g^\WM_{jk}}{\partial x^i} \Big)v^k dx^i\wedge dx^j + \sum_{i,j} g^\WM_{ij}dx^i \wedge dv^j.
  \end{equation}
  Due to the Hessian structure \eqref{eq:Lagriangian-two-form-on-TW:g-as-second-derivative}
  \begin{equation}
    \frac{\partial g^\WM_{ik}}{\partial x^j} = \frac{\partial^3 \what{\vphi}_\mrm{acc}}{\partial x^j \partial x^i \partial x^k} = \frac{\partial^3 \what{\vphi}_\mrm{acc}}{\partial x^i \partial x^j \partial x^k} = \frac{\partial g^\WM_{jk}}{\partial x^i}
  \end{equation}
  holds and the first sum in \eqref{eq:Lagriangian-two-form-on-TW:local-expr-Riemannian-theta_L} vanishes, resulting in the simplified expression
  \begin{equation}\label{eq:Lagriangian-two-form-on-TW:EQa}
    \om_L = \sum_{ij}g^\WM_{ij}dx^i \wedge dv^j.
  \end{equation}
  Suppose $A = (A', A''), B=(B', B'') \in T_{(W, V)} T\WM = \TW \times \TW$ with coordinates
  \begin{equation}
    A = \sum_i A^{'i} \frac{\partial }{\partial x^i} + \sum_i A^{''i} \frac{\partial }{\partial v^i} 
    \quad\text{and}\quad
    B = \sum_i B^{'i} \frac{\partial }{\partial x^i} + \sum_i B^{''i} \frac{\partial }{\partial v^i}.
  \end{equation}
  Evaluating the Lagriangian two-form \eqref{eq:Lagriangian-two-form-on-TW:EQa} we finally obtain
  \begin{equation}
    \om_L(A, B) = \sum_{ij} g^\WM_{ij}A^{'i}B^{''j} 
      - \sum_{ji}g^\WM_{ji}A^{''j}B^{'i} = g(A', B'') - g(A'', B').\qedhere
  \end{equation}
\end{proof}

Now that we have an explicit expression for the Lagrangian two-form $\om_L$, we are in a position to calculate an explicit representation of the Lagrangian vector field $X_E$.

\begin{proposition}\label{prop:Lagriangian-vf-on-W}
  The Lagrangian vector field $X_E$ on $T\WM$ associated to the Lagrangian \eqref{eq:L-by-G} at a point $(W, V) \in T\WM = \WM \times \TW$ is given by
  \begin{equation}\label{eq:Can-Sympl-Struct:Formulations-Lagrangian-Dynamics:Lvf-on-WM}
    X_E(W, V) = \bpm V \\ \frac{1}{2} \mc{A}(W, V) - \ggradg G(W) \epm.
  \end{equation}
\end{proposition}
\begin{proof}
  We directly use the definition \eqref{eq:def-Lagrangian-vector-field} of the Lagrangian vector field $X_E$. For this, let $B = (B', B'') \in T_{(W, V)} T\WM = \TW \times \TW$ be arbitrary and assume, $\gamma(t) = (W(t), V(t))$ is a smooth curve in $T\WM = \WM \times \TW$ with
  \begin{equation}
    \gamma(0) = (W(0), V(0)) = (W, V)\quad \text{and}\quad \dot{\gamma}(0) = (\dot{W}(0), \dot{V}(0)) = B = (B', B'').
  \end{equation}
  The time derivative of the potential $G$ is expressed via the Rimannian gradient
  \begin{equation}
    \tfrac{d}{dt} G(W(t))\big|_{t=0} = dG|_W[B'] = g_{W}\big(\ggradg G(W), B'\big).
  \end{equation}
  By \eqref{eq:Mechanics_of_Assignment_Flows:covariant-derivative-on-W}, the covariant derivative of $V(t)$ at $t=0$ is $D_t^gV(0) = B'' - \frac{1}{2}\mc{A}(W, V)$, resulting in
  \begin{subequations}
  \begin{align}
    \tfrac{d}{dt}\tfrac{1}{2}\big\|V(t)\big\|_{g}^2\big|_{t=0} &= \tfrac{1}{2} \tfrac{d}{dt} g_{W(t)}\big(V(t), V(t)\big)\big|_{t=0} = g_{W(0)}\big(V(0), D_t^g V(0)\big) \\
    &= g_{W}\big(V, B'' - \tfrac{1}{2}\mc{A}(W, V)\big).
  \end{align}
  \end{subequations}
  Putting everything together we obtain the following relation for the differential of the energy $E$ from \eqref{eq:Lagrangian-Systems:L-Riemann-mfld-Action-and-Energy}
  \begin{subequations}
  \begin{align}
    dE|_{(W, V)}[B] &= \tfrac{d}{dt} E(W(t), V(t))\big|_{t=0} = \tfrac{d}{dt}\tfrac{1}{2}\big\|V(t)\big\|_{g}^2\big|_{t=0} + \tfrac{d}{dt} G(W(t))\big|_{t=0}\\
    &= g_W\big(V, B''\big) - g_W\big(\tfrac{1}{2} \mc{A}(W, V) - \ggradg G(W), B'\big).\label{eq:Lagriangian-vf-on-W:calculation-dE}
  \end{align}
  \end{subequations}
  Writing $X_E = (X_E', X_E'') \in \TW \times \TW$ and comparing \eqref{eq:Lagriangian-vf-on-W:calculation-dE} with the above expression for $\om_L$ from Lemma~\ref{lem:Lagriangian-two-form-on-TW} shows $X_E'(W, V) = V$ and $X_E''(W, V) = \tfrac{1}{2} \mc{A}(W, V) - \ggradg G(W)$.
\end{proof}

Any solution curve $\gamma(t) = (W(t), V(t)) \in T\WM = \TW\times \TW$ of the Lagrangian dynamics induced by the Lagrangian vector field $X_E$ associated to the Lagrangian \eqref{eq:L-by-G} of Theorem~\ref{thm:Introduction:main-theorem} fulfills the ODE
\begin{equation}
  \bpm \dot{W}\\ \dot{V} \epm = X_E(W, V) = \bpm V \\ \frac{1}{2} \mc{A}(W, V) - \ggradg G(W) \epm.
\end{equation}
This form of the Hamiltonian ODE simply reflects the fact that this first-order dynamics on $T\WM$ is induced by a second-order ODE on $\WM$. Indeed, substituting $V = \dot{W}$ in the second component of $X_E$ results in
\begin{equation}
  \ddot{W} = \dot{V} = \tfrac{1}{2} \mc{A}(W, V) - \ggradg G(W) \quad \overset{\eqref{eq:Mechanics_of_Assignment_Flows:acceleration-on-W}}{\Leftrightarrow} \quad D^g_t \dot{W} = - \ggradg G(W),
\end{equation}
which we have already known to be satisfied for the base curve $W(t)$ by \eqref{eq:Introduction:main-theorem:E-L-Eq} of Theorem~\ref{thm:Introduction:main-theorem}.

\section{Conclusion}
\label{sec:Conclusion}

In this work, we generalized a previous result of uncoupled replicator equations from \cite{raju2018variational} to the case of coupled replicator equations. The viewpoint of Lagrangian mechanics on manifolds resulted in an interpretable Euler-Lagrange equation \eqref{eq:Introduction:main-theorem:E-L-Eq} and provided the mathematical tools to derive condition \eqref{eq:Introduction:main-theorem:characterization} for characterizing those affinity maps $F$ that result in  
critical points of the action functional \eqref{eq:Introduction:var-functional-L-on-W}. Accordingly, a constructed counterexample in terms of the specific affinity map \eqref{eq:counterexample-def-F} highlights that not all affinity maps $F$ lead to critical points.

The geometric mechanics perspective enabled the insight that, ignoring a set of starting points of measure zero, solutions to the assignment flow are  reparametrized geodesics of the Jacobi metric \eqref{eq:energy-0-Jacobi-metric}. Thus, in a certain sense, these solutions locally connect assignment states in an optimal way by realizing a shortest path. Finally, using the Legendre transformation, we calculated an explicit expression for the associated Hamiltonian system in terms of the corresponding Lagrangian system \eqref{eq:Can-Sympl-Struct:Formulations-Lagrangian-Dynamics:Lvf-on-WM}.
 
Our results provide a basis for exploring analogies to mathematical representations of interacting particle systems in theoretical physics in future work. This may further enhance our understanding of dynamical and learning systems that reveal structures in metric data.

\subsection*{Acknowledgements} 
This work is supported by the Deutsche Forschungsgemeinschaft (DFG, German Research Foundation) under Germany's Excellence Strategy EXC 2181/1 - 390900948 (the Heidelberg STRUCTURES Excellence Cluster) and the Transregional Colloborative Research Center CRC/TRR 191 (281071066).

\bibliographystyle{amsalpha}
\bibliography{journal_paper_mechanics}

\newcommand{\etalchar}[1]{$^{#1}$}
\providecommand{\bysame}{\leavevmode\hbox to3em{\hrulefill}\thinspace}
\providecommand{\MR}{\relax\ifhmode\unskip\space\fi MR }
\providecommand{\MRhref}[2]{%
  \href{http://www.ams.org/mathscinet-getitem?mr=#1}{#2}
}
\providecommand{\href}[2]{#2}
\begin{thebibliography}{MBP{\etalchar{+}}21}

\bibitem[Abb13]{abbo2013}
A.~Abbondandolo, \emph{Lectures on the free period {L}agrangian action
  functional}, J. Fixed Point Theory Appl. \textbf{13} (2013), no.~2, 397--430.

\bibitem[AC10]{amari2010information}
S.-I. Amari and A.~Cichocki, \emph{{Information Geometry of Divergence
  Functions}}, {Bulletin of the Polish Academy of Sciences. Technical Sciences}
  \textbf{58} (2010), no.~1, 183--195.

\bibitem[AJLS17]{ay2017information}
N.~Ay, J.~Jost, H.~V. L{\^e}, and L.~Schwachh{\"o}fer, \emph{{Information
  Geometry}}, Springer, 2017.

\bibitem[AM87]{abraham1987foundations}
R.~Abraham and J.~E. Marsden, \emph{{Foundations of Mechanics}}, 2nd ed.,
  Addison-Wesley Publ. Comp., Inc., 1987.

\bibitem[AN00]{amari2007methods}
S.-I. Amari and H.~Nagaoka, \emph{{Methods of Information Geometry}}, Amer.
  Math. Soc. and Oxford Univ. Press, 2000.

\bibitem[APSS17]{Astroem2017}
F.~Astr\"om, S.~Petra, B.~Schmitzer, and C.~Schn\"orr, \emph{{I}mage {L}abeling
  by {A}ssignment}, J.~Math.~Imag.~Vision \textbf{58} (2017), no.~2, 211--238.

\bibitem[CRBD18]{Chen:2018ab}
R.~T.~Q. Chen, Y.~Rubanova, J.~Bettencourt, and D.~Duvenaud, \emph{{Neural
  Ordinary Differential Equations}}, Proc. NeurIPS, 2018.

\bibitem[HR17]{Haber:2017aa}
E.~Haber and L.~Ruthotto, \emph{{Stable Architectures for Deep Neural
  Networks}}, Inverse Problems \textbf{34} (2017), no.~1, 014004.

\bibitem[HS03]{hofbauer2003}
J.~Hofbauer and K.~Sigmund, \emph{{Evolutionary Game Dynamics}}, Bull. Amer.
  Math. Soc. \textbf{40} (2003), no.~4, 479--519.

\bibitem[KAH{\etalchar{+}}15]{Kappes:2015aa}
J.H. Kappes, B.~Andres, F.A. Hamprecht, C.~Schn{\"o}rr, S.~Nowozin, D.~Batra,
  S.~Kim, B.X. Kausler, T.~Kr{\"o}ger, J.~Lellmann, N.~Komodakis,
  B.~Savchynskyy, and C.~Rother, \emph{{A Comparative Study of Modern Inference
  Techniques for Structured Discrete Energy Minimization Problems}}, Int. J.
  Computer Vision \textbf{115} (2015), no.~2, 155--184.

\bibitem[Lee13]{lee2013smooth}
John~M Lee, \emph{{Smooth Manifolds}}, Springer, 2013.

\bibitem[MBP{\etalchar{+}}21]{Minaee:2021tx}
S.~Minaee, Y.~Y. Boykov, F.~Porikli, A.~J. Plaza, N.~Kehtarnavaz, and
  D.~Terzopoulos, \emph{{Image Segmentation Using Deep Learning: A Survey}},
  IEEE Trans. Pattern Anal. Mach. Intell.,
  https://ieeexplore.ieee.org/document/9356353 (2021).

\bibitem[Mn97]{mane1997}
R.~Ma\~{n}\'{e}, \emph{Lagrangian flows: the dynamics of globally minimizing
  orbits}, Bol. Soc. Brasil. Mat. (N.S.) \textbf{28} (1997), no.~2, 141--153.

\bibitem[RK18]{raju2018variational}
V.~Raju and P.~S. Krishnaprasad, \emph{{A Variational Problem on the
  Probability Simplex}}, IEEE Conf. on Decision and Control (CDC), 2018,
  pp.~3522--3528.

\bibitem[San10]{sandholm2010}
W.~H. Sandholm, \emph{{Population Games and Evolutionary Dynamics}}, MIT press,
  2010.

\bibitem[Sch20]{Schnorr2019aa}
C.~Schn\"{o}rr, \emph{{Assignment Flows}}, {Handbook of Variational Methods for
  Nonlinear Geometric Data} (P.~Grohs, M.~Holler, and A.~Weinmann, eds.),
  Springer, 2020, pp.~235---260.

\bibitem[SS21]{Savarino2019ab}
F.~Savarino and C.~Schn\"{o}rr, \emph{{Continuous-Domain Assignment Flows}},
  Europ. J. Appl. Math. \textbf{32} (2021), no.~3, 570--597.

\bibitem[SY97]{shima1997geometry}
H.~Shima and K.~Yagi, \emph{{Geometry of Hessian Manifolds}}, Differential
  Geometry and its Applications \textbf{7} (1997), no.~3, 277--290.

\bibitem[ZSPS20]{Zeilmann:2020aa}
A.~Zeilmann, F.~Savarino, S.~Petra, and C.~Schn{\"{o}}rr, \emph{{Geometric
  Numerical Integration of the Assignment Flow}}, Inverse Problems \textbf{36}
  (2020), no.~3, 034004 (33pp).

\bibitem[ZZS21]{Zern:2020aa}
A.~Zern, A.~Zeilmann, and C.~Schn\"{o}rr, \emph{{Assignment Flows for Data
  Labeling on Graphs: Convergence and Stability}}, Information Geometry (in
  press; arXiv:2002.11571, 2021).

\end{thebibliography}

\end{document}